\documentclass[11pt, reqno]{amsart}

\usepackage[T2A, T1]{fontenc}    
\usepackage[utf8]{inputenc}     
\usepackage[russian,english]{babel}

\usepackage[margin=.9in, a4paper, foot=.3in]{geometry}

\makeatletter
\def\section{\@startsection{section}{1}%
  \z@{-.7\linespacing\@plus -\linespacing}{.5\linespacing}%
  {\normalfont\scshape\centering}}
\def\subsection{\@startsection{subsection}{2}%
  \z@{-.5\linespacing\@plus -.7\linespacing}{.5em}%
  {\normalfont\bfseries\mathversion{bold}}}
\makeatother

\usepackage{graphicx, xcolor}

\definecolor{darkblue}{rgb}{0,0,.5}
\definecolor{darkred}{rgb}{.5,0,0}
\definecolor{darkgreen}{rgb}{0,0.5,0}

\usepackage{hyperref}
\hypersetup {
    pdfstartview=FitH,
    colorlinks,
    citecolor=darkgreen,
    linkcolor=darkred,
    urlcolor=darkblue
}

\usepackage[noBBpl]{mathpazo}
\usepackage[mathcal]{euscript}
\usepackage{amssymb}
\usepackage{bm, bbm, mathtools}
\allowdisplaybreaks
\numberwithin{equation}{section}

\usepackage{cite}

\newcommand {\bil}[2]{(#1|#2)}

\newcommand {\ldbr}{[\![}
\newcommand {\oo}{{\overline 1}}
\newcommand {\oz}{{\overline 0}}
\newcommand {\rdbr}{]\!]}
\newcommand {\str}{\mathrm{str}}

\newcommand {\mbar}[3]{\hskip #2 \overline{\hskip -#2 #1 \hskip -#3} \hskip #3}
\newcommand {\obmE}{\mbar{\bm E}{.15em}{.05em}}

\newcommand {\rme}{\mathrm e}

\newcommand {\bbC}{\mathbb C}
\newcommand {\bbH}{\mathbb H}
\newcommand {\bbK}{\mathbb K}
\newcommand {\bbM}{\mathbb M}
\newcommand {\bbX}{\mathbb X}
\newcommand {\bbZ}{\mathbb Z}

\newcommand {\calE}{\mathcal E}
\newcommand {\calK}{\mathcal K}
\newcommand {\calR}{\mathcal R}

\newcommand {\gothg}{\mathfrak g}
\newcommand {\gothh}{\mathfrak h}
\newcommand {\gothk}{\mathfrak k}
\newcommand {\hgothh}{\widehat{\mathfrak h}}

\newcommand {\glto}{\mathfrak{gl}_{2|1}}
\newcommand {\hlslto}{\widehat{\mathcal{L}}(\mathfrak{sl}_{2|1})}
\newcommand {\lslm}{\mathcal L(\mathfrak{sl}_M)}
\newcommand {\lslto}{\mathcal{L}(\mathfrak{sl}_{2|1})}
\newcommand {\slm}{\mathfrak{sl}_M}
\newcommand {\slto}{\mathfrak{sl}_{2|1}}
\newcommand {\tlslto}{\widetilde{\mathcal{L}}(\mathfrak{sl}_{2|1})}
\newcommand {\uqg}{\mathrm U_q(\mathfrak g)}
\newcommand {\uqglmn}{\mathrm U_q(\mathfrak{gl}_{M|N})}
\newcommand {\uqglm}{\mathrm U_q(\mathfrak{gl}_M)}
\newcommand {\uqglto}{\mathrm U_q(\mathfrak{gl}_{2|1})}
\newcommand {\uqglthree}{\mathrm U_q(\mathfrak{gl}_3)}
\newcommand {\uqgltwo}{\mathrm U_q(\mathfrak{gl}_2)}
\newcommand {\uqslto}{\mathrm U_q(\mathfrak{sl}_{2|1})}
\newcommand {\uqlslm}{\mathrm U_q(\mathcal L(\mathfrak{sl}_M))}
\newcommand {\uqlslmn}{\mathrm U_q(\mathcal L(\mathfrak{sl}_{M|N}))}
\newcommand {\uqlslto}{\mathrm U_q(\mathcal L(\mathfrak{sl}_{2|1}))}

\newcommand {\uqlslthree}{\mathrm U_q(\mathcal L(\mathfrak{sl}_3))}
\newcommand {\uqlsltwo}{\mathrm U_q(\mathcal L(\mathfrak{sl}_2))}

\DeclareMathOperator {\End}{End}
\DeclareMathOperator {\id}{id}

\title{Basic monodromy operator for quantum superalgebra}

\author{A. V. Razumov}
\address{Institute for High Energy Physics, NRC ``Kurchatov Institute", 142281 Protvino, Mos\-cow region, Russia}

\email{Alexander.Razumov@ihep.ru}

\begin{document}

\addtolength {\jot}{3pt}

\begin{abstract}
We derive the explicit form of the basic monodromy operator for the quantum loop superalgebra $\mathrm{U}_q(\mathcal{L}(\mathfrak{sl}_{2|1}))$. Two significant additional results emerge from this derivation: simple expressions for the generating functions of the the images of the root vectors of $\mathrm{U}_q(\mathcal{L}(\mathfrak{sl}_{2|1}))$ under the Drinfeld homomorphism and explicit expressions for certain central elements of the quantum superalgebra $\mathrm{U}_q(\mathfrak{gl}_{2|1})$. Furthermore, we establish the relationship between these central elements and those obtained by using the Drinfeld partial trace method.
\end{abstract}

\maketitle


\section{Introduction}

Quantum algebras, also known as quantum groups, play a significant role in the study of quantum integrable systems. The general notion of a quantum algebra $\uqg$ was introduced by Drinfeld and Jimbo \cite{Dri85, Jim85} for the case when $\gothg$ is a Kac--Moody algebra with a symmetrizable generalized Cartan matrix. In the quantum algebraic approach, all the objects that describe the model and are related to its integrability originate from the universal $R$-matrix of the corresponding quantum algebra. We refer to these objects as integrability objects.

The universal $R$-matrix is an element of the tensor product of two copies of the quantum algebra under consideration. The integrability objects are constructed by selecting representations of the quantum algebra for the factors of that tensor product. When the same representation of the quantum algebra is chosen for both factors, we obtain an $R$-matrix. By using different representations for the factors, a more general integrable object, known as a monodromy operator, is produced. An even broader class of integrable objects can be constructed by applying a homomorphism from the quantum algebra to a certain associative algebra for one of the tensor product factors. For examples, see the papers \cite{BazLukZam96, BazLukZam97, BazLukZam99, BazHibKho02, BooGoeKluNirRaz10, BooGoeKluNirRaz11, BooGoeKluNirRaz13, Raz13, BooGoeKluNirRaz14a}.

Consider, for example, the case of integrable systems associated with the quantum algebra $\uqlslm$, where $\lslm$ denotes the loop algebra of the Lie algebra $\slm$. As discovered by Jimbo \cite{Jim86a}, there exists an algebra homomorphism $\epsilon$ from the quantum algebra $\uqlslm$ to the quantum algebra $\uqglm$. Using this homomorphism, for any given representation $\pi$ of the quantum algebra $\uqglm$, we define the representation $\pi \circ \epsilon$ of the quantum algebra $\uqlslm$. This type of representation is known as evaluation representations.

A monodromy operator, known as the basic monodromy operator, is defined by applying the homomorphism $\epsilon$ to the first factor of the universal $R$-matrix and a fixed representation of $\uqlslm$ to the second. This construction allows us to investigate the integrable systems related to all evaluation representations simultaneously. 

The explicit derivations of the basic monodromy operators for the quantum algebras $\uqlsltwo$ and $\uqlslthree$ are presented in the papers \cite{Raz13} and \cite{NirRaz16a}. Two significant results emerge from these derivations: the explicit forms of certain central elements of the quantum algebras $\uqgltwo$ and $\uqglthree$, and simple expressions for the generating functions of the images of the root vectors of $\uqlsltwo$ and $\uqlslthree$ under the Jimbo homomorphism. These generating functions are crucial for deriving the $\ell$-weights of the evaluation representations \cite{BooGoeKluNirRaz16, BooGoeKluNirRaz17b}. Along with the $\ell$-weights of the corresponding $q$-oscillator representations \cite{NirRaz17a, NirRaz17b}, they are used to prove the functional relations for the systems under consideration \cite{Raz21, Raz21a}.

By appropriately generalizing the defining relations of quantum algebras, one can obtain quantum superalgebras associated with Lie superalgebras \cite{Yam91, Yam94}. Numerous papers have been devoted to investigating the corresponding integrable systems (see, for example, the paper \cite{BazTsu08} and the references therein). However, a quantum algebraic derivation of the explicit form of the basic monodromy operator for such systems has yet to be provided. It is worth noting that in the paper \cite{Zha92}, the basic monodromy operator for the quantum superalgebra $\uqlslmn$ was constructed via the Baxterization procedure from the monodromy operators of the quantum algebra $\uqglmn$. In this paper, we present a direct quantum algebraic derivation of the basic monodromy operator for the quantum superalgebra $\uqlslto$. As with $\uqlslthree$, we also obtain two additional important results: the explicit forms of certain central elements and simple expressions for the images of the generating functions of the root vectors of the quantum superalgebra $\uqlslto$ under the Jimbo homomorphism.

The structure of the paper is as follows. In Subsection \ref{s:2.1} we review the necessary definitions and essential facts related to the Lie superalgebras $\glto$ and $\slto$. In Subsection \ref{s:2.2} we introduce the quantum superalgebra $\uqglto$ and briefly discuss its representations. Subsection \ref{s:3.1} examines the affine Lie superalgebra $\hlslto$, followed by the definition of the quantum loop superalgebra $\uqlslto$ in Subsection \ref{s:3.2}. In Subsection \ref{s:3.3}, we present the Jimbo homomorphism from $\uqlslto$ to $\uqglto$ and discuss the modules and representations of $\uqslto$ generated through this homomorphism. We then describe the expression for the universal $R$-matrix proposed by Khoroshkin and Tolstoy \cite{KhoTol91, KhoTol93a, KhoTol94a}. To this end, in Subsection \ref{s:3.4} we give the relations used to construct the Cartan--Weyl generators of $\uqlslto$, and in Subsection \ref{s:3.5} provide the explicit expression for the universal $R$-matrix. The general definition of an integrability object is given in Subsection \ref{s:4.1}. Subsection \ref{s:4.2} introduces the notion of a monodromy operator and the basic monodromy operator. Subsection \ref{s:4.3} discusses the structure of the basic monodromy operator, derived from the structure of the universal $R$-operator. In Subsection \ref{s:4.4}, we propose an ansatz for the basic monodromy operator and obtain the explicit expression for this operator. In Section \ref{s:5} we establish the relationship between the central elements introduced in Subsection \ref{s:4.4} and those obtained by using the Drinfeld partial trace method.

We use the notation
\begin{equation*}
\kappa_q = q - q^{-1},
\end{equation*}
so that the $q$-deformation of the integer $n$ is
\begin{equation*}
[n]_q = \frac{q^n - q^{-n}}{q - q^{-1}} = \kappa_q^{-1} (q^n - q^{-n}).
\end{equation*}
We define the $q$-exponential by the equation
\begin{equation*}
\exp_q(x) = \sum_{n = 0}^\infty \frac{x^n}{(n)_q!} ,
\end{equation*}
where
\begin{equation*}
(n)_q! = (n)_q (n - 1)_q \ldots (1)_q,\qquad (n)_q = \frac{q^n - 1 }{q - 1}.
\end{equation*}
For calculations in $\uqglto$ we use the Wolfram Mathematica system.

\section{\texorpdfstring{Quantum superalgebra $\uqglto$}{Quantum superalgebra Uq(gl2|1)}} 

\subsection{\texorpdfstring{Lie superalgebras $\glto$ and $\slto$}{Lie superalgebras gl2|1 and sl2|1}} \label{s:2.1}

Let $V$ be a $\bbZ_2$-graded vector space such that $\dim V = (2,1)$. For convenience, we denote the Lie superalgebra $\mathfrak{gl}(V)$ as $\glto$.\footnote{The necessary information on the superalgebras and Lie superalgebras can be found in Appendix \ref{a:a}.} Let $(\mathbbm v_i)_{i = 1,2,3}$ be a basis of $V$ consisting of homogeneous elements, with
\begin{equation*}
[\mathbbm v_1] = [\mathbbm v_2] = \oz, \qquad [\mathbbm v_3] = \oo.
\end{equation*}
It is convenient to use the notation
\begin{equation*}
[i] = [\mathbbm{v}_i].
\end{equation*}
The elements $\bbM_{i j} \in \End(V)$ for $i, j = 1, 2, 3$, defined by the equation
\begin{equation*}
\bbM_{i j} \mathbbm{v}_k = \mathbbm{v}_i \delta_{j k},
\end{equation*}
form a basis of $\End(V)$ consisting of homogeneous elements. We have
\begin{equation*}
\bbM_{i j} \bbM_{k l} = \delta_{j k} \bbM_{i l}.
\end{equation*}
It is also clear that
\begin{equation*}
[\bbM_{i j}] = [i] + [j].
\end{equation*}

Denote by $\gothk$ the Cartan subalgebra of the Lie superalgebra $\glto$, generated by the elements $\bbK_i = \bbM_{i i}$ for $i = 1, 2, 3$, which form its standard basis.

Let the elements $\Xi_i$ for $i = 1, 2, 3$, form the dual basis of $\gothk^*$. The elements
\begin{equation*}
\alpha_{i j} = \Xi_i -\Xi_j, \qquad i \ne j,
\end{equation*}
form the root system $\Delta$ of $\glto$. We choose the elements
\begin{equation*}
\alpha_i = \alpha_{i, \, i + 1}, \qquad i = 1, 2,
\end{equation*}
as the simple root. In this case the system of positive roots $\Delta_+$ is formed by the roots
\begin{equation*}
\alpha_{i j} =  \sum_{k = 1}^{j - 1} \alpha_k, \qquad 1 \le i < j \le 3.
\end{equation*}
Note that the root $\alpha_1$ is even and $\alpha_2$ is odd, see the paper \cite{Kac77} for the definition. Denoting the parity of a root $\gamma \in \Delta$ by $[\gamma]$, we see that
\begin{equation*}
[\alpha_{i j}] = [\bbM_{i j}] = [i] + [j].
\end{equation*}
In fact, we consider the parity of positive roots as an additive mapping from $\Delta_+$ to $\bbZ_2$.

We define a symmetric bilinear form $\bil {\cdot}{\cdot}$ on $\gothk$ by the equation
\begin{equation*}
\bil {\bbK_i}{\bbK_j} = (-1)^{[i]} \delta_{i j}.
\end{equation*}
The form $\bil {\cdot}{\cdot}$ is nondegenerate and induces a symmetric bilinear form on $\gothk^*$. It is evident that
\begin{equation*}
\bil {\Xi_i}{\Xi_j} = (-1)^{[i]} \delta_{i j}.
\end{equation*}

The supertrace of an element $\bbX = \sum_{i, j = 1}^3 X_{i j} \, \bbM_{i j} \in \glto$ is defined as
\begin{equation*}
\str \, \bbX =  \sum_{i = 1}^3 (-1)^{[i]} X_{i i}.
\end{equation*}
The special linear Lie superalgebra $\slto$ is a Lie subsuperalgebra of $\glto$ formed by the elements with zero supertrace. As a Cartan subalgebra $\gothh$ of $\slto$, we take the subalgebra of $\gothk$, which consists of all elements  with zero supertrace. The standard basis of $\gothh$ is formed by the elements $\bbH_i$ for $i = 1, 2$, defined as
\begin{equation*}
\bbH_1 = \bbK_1 - \bbK_2, \qquad \bbH_2 = \bbK_2 + \bbK_3.
\end{equation*}
As the positive and negative roots of $\slto$ we take the restriction of the roots of $\glto$ to $\gothh$. The numbers
\begin{equation*}
a_{i j} = \alpha_j(H_i), \qquad i, j = 1, 2,
\end{equation*}
are the entries of the Cartan matrix $A$ of the Lie superalgebra $\slto$. We have
\begin{equation*}
A = (a_{i j}) = \left( \begin{array}{rr}
2 & -1 \\
-1 & 0
\end{array} \right).
\end{equation*}
The restriction of the bilinear form $\bil {\cdot}{\cdot}$ to $\gothh$ is nondegenerate and symmetric. We have
\begin{equation*}
\bil {\bbH_i}{\bbH_j} =  a^{}_{i j}.
\end{equation*}
and for the induced bilinear form on $\gothh^*$
\begin{equation*}
\bil {\alpha_i}{\alpha_j} = a_{i j}.
\end{equation*}
One can show that
\begin{align*}
& \bil{\alpha_{1 2}}{\alpha_{1 2}} = 2, && \bil{\alpha_{1 2}}{\alpha_{2 3}} = -1, && \bil{\alpha_{1 2}}{\alpha_{1 3}} = 1, \\
& \bil{\alpha_{2 3}}{\alpha_{1 2}} = -1, && \bil{\alpha_{2 3}}{\alpha_{2 3}} = 0, && \bil{\alpha_{2 3}}{\alpha_{1 3}} = -1, \\
& \bil{\alpha_{1 3}}{\alpha_{1 2}} = 1, && \bil{\alpha_{1 3}}{\alpha_{2 3}} = -1, && \bil{\alpha_{1 3}}{\alpha_{1 3}} = 0.
\end{align*}

\subsection{\texorpdfstring{Definition of quantum superalgebra $\uqglto$}{Definition of quantum superalgebra Uq(gl2|1}} \label{s:2.2}

Let $\hbar$ be an indeterminate and $q = \exp(\hbar)$. We define the quantum superalgebra $\uqglto$ as a unital associative $\bbC[[\hbar]]$-superalgebra generated by the elements\footnote{We use bold letters to denote the generators of quantum algebras and  objects related to them.}
\begin{equation*}
\bm K_i, \quad i = 1, 2, 3, \qquad \bm E_i, \quad \bm F_i, \quad i = 1, 2,
\end{equation*}
which satisfy the corresponding defining relations. The $\bbZ_2$-grading of $\uqglto$ is defined on the generators as follows:
\begin{equation*}
[\bm K_1] = [\bm K_2] = [\bm K_3] = \overline 0, \qquad [\bm E_1] = [\bm F_1] = \oz, \qquad [\bm E_2] = [\bm F_2] = \oo.
\end{equation*}

Before presenting the defining relations of $\uqglto$, we introduce the notion of the $q$-supercom\-mutator. The abelian group
\begin{equation*}
Q = \bigoplus_{i = 1}^2 \bbZ \, \alpha_i.
\end{equation*}
is called the root lattice of $\glto$. We denote by $Q_+$ and $Q_-$ the sublattices
\begin{equation*}
Q_+ = \bigoplus_{i = 1}^2 \bbZ_{\ge 0} \, \alpha_i, \qquad Q_- = \bigoplus_{i = 1}^2 \bbZ_{\le 0} \, \alpha_i.
\end{equation*}
Assuming that
\begin{equation*}
\bm K_i \in \uqglto_0, \qquad \bm E_i \in \uqglto_{\alpha_i}, \qquad \bm F_i \in \uqglto_{-\alpha_i},
\end{equation*}
we endow $\uqglto$ with a $Q$-grading.

For any $a \in \uqglto_\alpha$ and $b \in \uqglto_\beta$, we define their $q$-supercommutator by the equation
\begin{gather*}
\ldbr a, \, b \rdbr =  a b - (-1)^{[a] [b]} q^{- \bil \alpha \beta} b a = a b - (-1)^{[\alpha] [\beta]} q^{- \bil \alpha \beta} b a
\intertext{if $\alpha, \beta \in Q_+$, by the equation}
\ldbr a, \, b \rdbr =  a b - (-1)^{[a] [b]} q^{\bil \alpha \beta} b a = a b - (-1)^{[\alpha] [\beta]} q^{\bil \alpha \beta} b a
\intertext{if $\alpha, \beta \in Q_-$, and by the equation}
\ldbr a, \, b \rdbr =  a b - (-1)^{[a] [b]} b a = a b - (-1)^{[\alpha] [\beta]} b a
\end{gather*}
if $\alpha \in Q_+$ and $\beta \in Q_-$, or $\alpha \in Q_-$ and $\beta \in Q_+$.

The defining relations of the quantum superalgebra $\uqglto$ have the form
\begin{gather*}
\ldbr \bm K_i, \, \bm K_j \rdbr = 0, \qquad i, j = 1, 2, 3, \\
\ldbr \bm K_i, \, \bm E_j \rdbr = \alpha_j(\bbK_i) \bm E_j, \qquad \ldbr \bm K_i, \, \bm F_j \rdbr = - \alpha_j(\bbK_i) \bm F_j, \qquad i = 1, 2, 3, \qquad j = 1, 2,  \\
\ldbr \bm E_i, \, \bm F_j \rdbr = \delta_{i j} \frac{q^{\bm H_i} - q^{- \bm H_i}}{q^{} - q^{- 1}}, \qquad i, j = 1, 2, \\
\ldbr \bm E_2, \, \bm E_2 \rdbr = 0, \qquad \ldbr \bm F_2, \, \bm F_2 \rdbr = 0, \\
\ldbr \bm E_1, \, \ldbr \bm E_1, \, \bm E_2 \rdbr \rdbr = 0, \qquad \ldbr \bm F_1, \, \ldbr \bm F_1, \, \bm F_2 \rdbr \rdbr = 0,
\end{gather*}
where
\begin{equation*}
\bm H_1 = \bm K_1 - \bm K_2, \qquad \bm H_2 = \bm K_2 + \bm K_3.
\end{equation*}

We define the elements
\begin{equation*}
\bm E_{1 3} = \ldbr \bm E_1, \, \bm E_2 \rdbr = \bm E_1 \bm E_2 - q \bm E_2 \bm E_1, \qquad \bm F_{1 3} = \ldbr \bm F_1, \, \bm F_2 \rdbr = \bm F_2 \bm F_1 - q^{-1} \bm F_1 \bm F_2
\end{equation*}
and, for consistency, denote
\begin{equation*}
\bm E_{1 2} = \bm E_1, \qquad \bm E_{2 3} = \bm E_2, \qquad \bm F_{1 2} = \bm F_1, \qquad \bm F_{2 3} = \bm F_2.
\end{equation*}
We observe that
\begin{equation*}
[\bm E_{i j}] = [i] + [j].
\end{equation*}
The root vectors $\bm E_{i j}$ and $\bm F_{i j}$  are linearly independent, and together with the elements $\bm K_i$, are referred to as the Cartan--Weyl generators of $\uqglto$. It appears that the ordered monomials constructed from the Cartan--Weyl generators form a Poincar\'e--Birkhoff--Witt basis of the quantum superalgebra $\uqglto$.

It is important that the quantum loop superalgebra $\uqglto$ is a Hopf superalgebra with respect to the comultiplication $\Delta$, the antipode $S$, and the counit $\varepsilon$ defined by the relations
\begin{gather*}
\Delta(\bm K_i) = \bm K_i \otimes 1 + 1 \otimes\bm K_i, \\
\Delta(\bm E_i) = \bm E_i \otimes 1 + q^{\bm H_i} \otimes \bm E_i, \qquad \Delta(\bm F_i) = \bm F_i \otimes q^{-\bm H_i} + 1 \otimes \bm F_i, \\
S(\bm K_i) = - \bm K_i, \qquad S(\bm E_i) = - q^{- \bm H_i} \bm E_i, \qquad S(\bm F_i) = - \bm F_i \, q^{\bm H_i}, \\
\varepsilon(\bm K_i) = 0, \qquad \varepsilon(\bm E_i) = 0, \qquad \varepsilon(\bm F_i) = 0.
\end{gather*}

\subsection{Modules and representations} \label{s:2.3}

We identify a triple $(\lambda_1, \, \lambda_2, \, \lambda_3)$ of complex numbers with an element $\lambda \in \gothk^*$ defined by the equation
\begin{equation*}
\lambda(\bbK_i) = \lambda_i,
\end{equation*}
and vice versa. A $\uqglto$-module $V$ is called a highest weight module of highest weight $\lambda \in \gothk^*$ if there exists a weight vector $v^\lambda \in V$ satisfying the relations
\begin{gather*}
\bm E_i v^\lambda = 0, \quad i = 1, 2, \qquad \bm K_i v^\lambda = \lambda_i v^\lambda, \quad i = 1, 2, 3, \\
\uqglto \, v^\lambda = V.
\end{gather*}
Given $\lambda \in \gothk^*$, let $\widetilde V^\lambda$ denote the corresponding Verma $\uqglto$-module. This is a highest weight module of highest weight $\lambda$. We denote by $\widetilde \pi^\lambda$ the representation of $\uqglto$ corresponding
to $\widetilde V^\lambda$. It is clear that $\widetilde V^\lambda$ and $\widetilde \pi^\lambda$ are infinite-dimensional. However, if the difference $\lambda_1 - \lambda_2$ is a non-negative integer, there is a maximal submodule such that the respective quotient module is finite-dimensional. We denote this quotient by $V^\lambda$ and the corresponding representation as $\pi^\lambda$.

In the present paper, we are especially interested in the finite-dimensional representation $\pi^{(1, \, 0, \, 0)}$ of $\uqglto$, which is realized on a $\bbZ_2$-graded vector space of $V$, introduced in Section \ref{s:2.1}, as follows:
\begin{gather*}
\pi^{(1,0,0)} (\bm K_1) = \bbK_1, \qquad \pi^{(1,0,0)} (\bm K_2) = \bbK_2, \qquad \pi^{(1,0,0)} (\bm K_3) = \bbK_3, \label{pib} \\
\pi^{(1,0,0)}(\bm E_{1 2}) =  \bbM_{1 2}, \quad \pi^{(1,0,0)}(\bm E_{2 3}) =  \bbM_{2 3}, \quad
\pi^{(1,0,0)}(\bm F_{1 2}) = \bbM_{2 1}, \quad \pi^{(1,0,0)}(\bm F_{2 3}) = \bbM_{3 2}. \label{pie}
\end{gather*}
It is straightforward to show that
\begin{equation*}
\pi^{(1,0,0)}(\bm E_{1 3}) =  \bbM_{1 3}, \qquad \pi^{(1,0,0)}(\bm F_{1 3}) = \bbM_{3 1}.
\end{equation*}
In what follows, we denote the representation $\pi^{(1,0,0)}$ simply as $\pi$. This representation is usually called the defining representation of $\uqglto$.

\section{\texorpdfstring{Quantum superalgebra $\uqlslto$}{Quantum superalgebra Uq(L(sl2|1))}}

\subsection{\texorpdfstring{On affine Lie superalgebra $\hlslto$}{On affine Lie superalgebra L(sl2|1)}} \label{s:3.1}

We begin with the necessary background on the nontwisted affine Lie superalgebra $\hlslto$. More information on affine Lie superalgebras can be found in the papers \cite{Leu86, Leu89}.

Let $\lslto$ denote the loop superalgebra of $\slto$, $\tlslto$ its standard central extension by a one-dimensional center $\bbC c$, and $\hlslto$ the Lie superalgebra obtained from $\tlslto$ by adding a natural derivation $d$. By definition,
\begin{equation*}
\hlslto = \lslto \oplus \bbC \, c \oplus \bbC \, d, 
\end{equation*}
and we use the space\footnote{Here we use the natural embedding of $\slto$ into $\hlslto$.} 
\begin{equation*}
\hgothh = \gothh \oplus \bbC \, c \oplus \bbC \, d
\end{equation*}
as the Cartan subsuperalgebra of $\hlslto$. We assume that $[c] = [d] = \overline 0$.

We denote
\begin{equation*}
h_1 = \bbH_1, \qquad h_2 = \bbH_2,
\end{equation*}
and define
\begin{equation*}
h_0 =  c - h_1 - h_2.
\end{equation*}
It is noteworthy that
\begin{equation*}
c = h_0 + h_1 + h_2.
\end{equation*}
We denote by $\delta$ the element of $\hgothh^*$ defined by the equations
\begin{equation*}
\delta(h_0) = 0, \qquad \delta(h_1) = 0, \qquad \delta(h_2) = 0, \qquad \delta(d) = 1,
\end{equation*}
and define the root $\alpha_0 \in \hgothh^*$ as
\begin{equation*}
\alpha_0 = \delta - \alpha_1 - \alpha_2.
\end{equation*}
We observe that
\begin{equation*}
\delta = \alpha_0 + \alpha_1 + \alpha_2.
\end{equation*}

The equation 
\begin{equation*}
a_{i j} = \alpha_j(h_i), \qquad i, j = 0, 1, 2,
\end{equation*}
gives the entries of the matrix
\begin{equation*}
A^{(1)} = (a_{i j})_{i, \, j = 0, \, 1, \, 2} = \left(\begin{array}{rrr}
0 & -1 & 1 \\
-1 & 2 & -1 \\
1 & -1 & 0
\end{array} \right),
\end{equation*}
which is the Cartan matrix of $\hlslto$.

We take $\alpha_i$, $i = 0, 1, 2$, as the simple roots. The full system of positive roots of the Lie superalgebra $\hlslto$ is
\begin{equation*}
\widehat \Delta_+ = \{\gamma + n \delta \mid  \gamma \in \Delta_+, \ n \in \bbZ_{\ge 0} \} \cup \{n \delta \mid n \in \bbZ_{>0} \} \cup \{(\delta - \gamma) + n \delta \mid  \gamma \in \Delta_+, \ n \in \bbZ_{\ge 0}\}.
\end{equation*}
The system of negative roots is $\widehat \Delta_- = - \widehat \Delta_+$, and the full root system is
\begin{equation*}
\widehat \Delta = \widehat \Delta_+ \cup \widehat \Delta_- 
= \{ \gamma + n \delta \mid \gamma \in \Delta, \ n \in \bbZ \} \cup \{n \delta \mid n \in \bbZ \setminus \{0\} \}.
\end{equation*}
The roots $n \delta$, $n \in \bbZ \setminus \{0\}$ are imaginary, all other roots are real.\footnote{A root $\gamma \in \Delta$ is called imaginary if $k \gamma \in \Delta$ for all $k \in \bbZ \setminus 0$, and real otherwise \cite{Leu86, Leu89}.} We assume that the root $\delta$ is even and extend the parity mapping from $\Delta_+$ to $\widehat \Delta_+$ by additivity. It is clear that
\begin{equation*}
[\gamma + n \delta] = [(\gamma - \delta) + n \delta] = [\gamma], \quad \gamma \in \Delta_+, \ n \in \bbZ_{\ge 0}, \qquad [n \delta] = \oz, \quad n \in \bbZ_{> 0} .
\end{equation*}

We fix a non-degenerate symmetric bilinear form $\bil{\cdot}{\cdot}$ on $\hgothh$ by the equations
\begin{equation*}
\bil{h_i}{h_j} = a^{}_{i j}, \qquad \bil{h_i}{d} = \delta^{}_{i 0}, \qquad \bil{d}{d} = 0, \qquad i, j = 0, 1, 2.
\end{equation*}
For the corresponding non-degenerate symmetric bilinear form on 
$\hgothh^*$ we have
\begin{equation*}
\bil{\alpha_i}{\alpha_j} = a_{i j}, \qquad i, j = 0, 1, 2.
\end{equation*}
It follows from this relation that
\begin{equation*}
\bil{\delta}{\gamma} = 0
\end{equation*}
for any $\gamma \in \widehat Q$.

\subsection{\texorpdfstring{Definition of $\uqlslto$}{Definition of Uq(L(sl2|1))}} \label{s:3.2}

Let again $\hbar$ be an indeterminate and $q = \exp (\hbar)$. The quantum loop superalgebra $\uqlslto$ is a unital associative $\bbC[[\hbar]]$-superalgebra generated by the elements
\begin{equation*}
\bm h_i, \quad \bm e_i, \quad \bm f_i, \quad i = 0, 1, 2.
\end{equation*}
The $\bbZ_2$-grading of $\uqlslto$ is defined on the generators as follows:
\begin{equation*}
[\bm h_0] = [\bm h_1] = [\bm h_2] = \overline 0, \qquad [\bm e_0] = [\bm f_0] = \oo, \qquad [\bm e_1] = [\bm f_1] = \oz, \qquad [\bm e_2] = [\bm f_2] = \oo.
\end{equation*}
The abelian group
\begin{equation*}
\widehat Q = \bigoplus_{i = 0}^2 \bbZ \, \alpha_i.
\end{equation*}
is called the root lattice of $\hlslto$. We denote by $\widehat Q_+$ and $\widehat Q_-$ the sublattices
\begin{equation*}
\widehat Q_+ = \bigoplus_{i = 0}^2 \bbZ_{\ge 0} \, \alpha_i, \qquad \widehat Q_- = \bigoplus_{i = 0}^2 \bbZ_{\le 0} \, \alpha_i.
\end{equation*}
Assuming that
\begin{equation}
\bm h_i \in \uqglto_0, \qquad \bm e_i \in \uqglto_{\alpha_i}, \qquad \bm f_i \in \uqglto_{-\alpha_i}, \qquad i = 0, 1, 2, \label{qxeifi}
\end{equation}
we endow $\uqlslto$ with a $\widehat Q$-grading. 

We define the $q$-supercommutator for the quantum loop superalgebra $\uqlslto$ in the same manner as for the quantum superalgebra $\uqglto$ in Section \ref{s:2.2}. The defining relations of the quantum loop superalgebra $\uqlslto$ look as
\begin{gather*}
\bm h_0 + \bm h_1 + \bm h_2 = 0, \\
\ldbr \bm h_i, \, \bm e_j \rdbr = a_{i j} \, \bm e_i, \qquad \ldbr \bm h_i, \, \bm f_j \rdbr = - a_{i j} \bm f_i, \qquad i, j = 0, 1, 2, \\
\ldbr \bm e_0, \, \bm e_0 \rdbr = 0, \qquad \ldbr \bm f_0, \, \bm f_0 \rdbr = 0, \qquad \ldbr \bm e_2, \, \bm e_2 \rdbr = 0, \qquad \ldbr \bm f_2, \, \bm f_2 \rdbr = 0,\\
\ldbr \bm e_i, \, \bm f_j \rdbr = \delta_{i j} \, \frac{q_i^{\bm h_i} - q_i^{- \bm h_i}}{q^{\mathstrut}_i - q_i^{-1}}, \qquad i, j = 0, 1, 2, \\
\ldbr \bm e_1, \, \ldbr \bm e_1, \, \bm e_0 \rdbr \rdbr = 0, \qquad \ldbr \bm f_1, \, \ldbr \bm f_1, \, \bm f_0 \rdbr \rdbr = 0, \\
\ldbr \bm e_1, \, \ldbr \bm e_1, \, \bm e_2 \rdbr \rdbr = 0, \qquad \ldbr \bm f_1, \, \ldbr \bm f_1, \, \bm f_2 \rdbr \rdbr = 0, \\
\ldbr \bm e_0, \, \ldbr \bm e_2, \, \ldbr \bm e_0, \, \ldbr \bm e_2, \, \bm e_1 \rdbr \rdbr \rdbr \rdbr = \ldbr \bm e_2, \, \ldbr \bm e_0, \, \ldbr \bm e_2, \, \ldbr \bm e_0, \, \bm e_1 \rdbr \rdbr \rdbr \rdbr, \\
\ldbr \bm f_0, \, \ldbr \bm f_2, \, \ldbr \bm f_0, \, \ldbr \bm f_2, \, \bm f_1 \rdbr \rdbr \rdbr \rdbr = \ldbr \bm f_2, \, \ldbr \bm f_0, \, \ldbr \bm f_2, \, \ldbr \bm f_0, \, \bm f_1 \rdbr \rdbr \rdbr \rdbr,
\end{gather*}
see, for example, the paper \cite{Yam99}.

The quantum loop superalgebra $\uqlslto$ is a Hopf superalgebra with respect to the comultiplication $\Delta$, the antipode $S$, and the counit $\varepsilon$ defined by the relations
\begin{gather*}
\Delta(\bm h_i) = \bm h_i \otimes 1 + 1 \otimes \bm h_i, \label{dqx} \\
\Delta(\bm e_i) = \bm e_i \otimes 1 + q^{\bm h_i} \otimes \bm e_i, \qquad \Delta(\bm f_i) = \bm f_i \otimes q^{-\bm h_i} + 1 \otimes \bm f_i, \\
S(\bm h_i) = - \bm h_i, \qquad S(\bm e_i) = - q^{- \bm h_i} \bm e_i, \qquad S(\bm f_i) = - \bm f_i \, q^{\bm h_i}, \\
\varepsilon(\bm h_i) = 0, \qquad \varepsilon(\bm e_i) = 0, \qquad \varepsilon(\bm f_i) = 0. 
\end{gather*}

\subsection{Modules and representations} \label{s:3.3}

To define representations of $\uqlslto$, we use the Jimbo homomorphism 
\begin{equation*}
\epsilon \colon \uqlslto \to \uqglto
\end{equation*}
which is defined by the relations
\begin{align}
& \epsilon(\bm h_0) = - \bm K_3 - \bm K_1, &&  \epsilon(\bm h_1) = \bm K_1 -\bm  K_2 && \epsilon(\bm h_2) = \bm K_2 + \bm K_3, \label{j2h} \\
& \epsilon(\bm e_0) = - \bm F_{1 3} \, q^{\bm K_1 - \bm K_3}, && 
\epsilon(\bm e_1) = \bm E_{1 2}, && \epsilon(\bm e_2) = \bm E_{2 3}, \label{j2e} \\
& \epsilon(\bm f_0) =  q^{- \bm K_1 + \bm K_3} \bm E_{1 3}, && 
\epsilon(\bm f_1) = \bm F_{1 2}, && \epsilon(\bm f_2) = \bm F_{2 3}. \label{j2f}
\end{align}

Using the Jimbo homomorphism and starting with the infinite-dimensional representation $\widetilde \pi^\lambda$ of $\uqglto$, we define the infinite-dimensional representation
\begin{equation*}
\widetilde \varphi^\lambda = \widetilde \pi^\lambda \circ \epsilon
\end{equation*}
of $\uqlslto$. Recall that when $\lambda_1 - \lambda_2$ is a non-negative integer, the representation $\widetilde \pi^\lambda$ has the maximal subrepresentation such that the corresponding quotient representation, denoted by $\pi^\lambda$, is finite-dimensional. We define the corresponding representation of $\uqlslto$ as
\begin{equation*}
\varphi^\lambda = \pi^\lambda \circ \epsilon.
\end{equation*}
With a slight abuse of notation, we denote the corresponding modules as $\widetilde V^\lambda$ and $V^\lambda$.

For the representation $\varphi^{(1, \, 0, \, 0)}$ we have
\begin{align}
& \varphi^{(1, \, 0, \, 0)}(\bm h_0) = - \bbK_3 - \bbK_1, &&  \varphi^{(1, \, 0, \, 0)}(\bm h_1) = \bbK_1 - \bbK_2 && \varphi^{(1, \, 0, \, 0)}(\bm h_2) = \bbK_2 + \bbK_3, \label{fh} \\
& \varphi^{(1, \, 0, \, 0)}(\bm e_0) = - q \, \bbM_{3 1}, && 
\varphi^{(1, \, 0, \, 0)}(\bm e_1) =  \bbM_{1 2}, && \varphi^{(1, \, 0, \, 0)}(\bm e_2) =  \bbM_{2 3}, \label{fe} \\
& \varphi^{(1, \, 0, \, 0)}(\bm f_0) = q^{-1}  \bbM_{1 3}, && 
\varphi^{(1, \, 0, \, 0)}(\bm f_1) = \bbM_{2 1}, && \varphi^{(1, \, 0, \, 0)}(\bm f_2) = \bbM_{3 2}. \label{ff}
\end{align}
In what follows, we denote the representation $\varphi^{(1, \, 0, \, 0)}$ and the module $V^{(1, \, 0, \, 0)}$ simply as $\varphi$ and $V$, respectively.

\subsection{Cartan--Weyl generators} \label{s:3.4}

Keeping in mind the $\widehat Q$-grading of $\uqlslto$ defined by equation~(\ref{qxeifi}), we say that an element $a$ of $\uqlslto$ is a root vector corresponding to a root $\gamma \in \widehat \Delta$ if $a \in \uqlslto_\gamma$. In particular, the generators $\bm e_i$ and $\bm f_i$ are root vectors corresponding to the roots $\alpha_i$ and $- \alpha_i$.

One can construct linearly independent root vectors corresponding to all roots from $\widehat \Delta \subset \widehat Q$. They are called the Cartan--Weyl generators of $\uqlslto$. If some ordering of roots is chosen, then appropriately ordered monomials constructed from these vectors and the generators $\bm h_i$ form a Poincar\'e--Birkhoff--Witt basis of $\uqlslto$. To construct the root vectors we follow the procedure, based on normal ordering of positive roots, proposed by Khoroshkin and Tolstoy \cite{KhoTol93a, KhoTol94, KhoTol94a, Raz23}.

For the case of a finite dimensional simple Lie superalgebra, an order relation $\prec$ is called a normal ordering of positive roots if each nonsimple positive root $\gamma = \alpha + \beta$, where $\alpha \prec \beta$, is positioned between $\alpha$ and $\beta$. In our case, we assume additionally that all imaginary roots follow each other in an arbitrary order, and that
\begin{equation}
\alpha + k \delta \prec m \delta \prec (\delta - \beta) + n \delta \label{akd}
\end{equation}
for any $\alpha, \beta \in \Delta_+$, $k, n \in \bbZ_{\ge 0}$ and $m \in \bbZ_{>0}$ \cite{KhoTol91, KhoTol93a, KhoTol94, KhoTol94a, Raz23}.

Assume that some normal ordering of positive roots is chosen.  We say that a pair $(\alpha, \, \beta)$ of positive roots generates a root $\gamma$ if $\gamma = \alpha + \beta$ and $\alpha \prec \beta$. A pair of positive roots $(\alpha, \, \beta)$ generating a root $\gamma$ is called minimal if there is no other pair of positive roots $(\alpha', \, \beta')$ generating $\gamma$ such that $\alpha \prec \alpha' \prec \beta' \prec \alpha$.

It is convenient to denote a root vector corresponding to a positive root $\gamma$ by $e_\gamma$, and a root vector corresponding to a negative root $- \gamma$ by $f_\gamma$. Following \cite{KhoTol93a, KhoTol94, KhoTol94a, Raz23}, we define root vectors using the following inductive procedure. Given a root $\gamma \in \widehat \Delta_+$, let $(\alpha, \, \beta)$ be a minimal pair of positive roots generating $\gamma$. Now, if the root vectors $e_\alpha$, $e_\beta$ and $f_\alpha$, $f_\beta$ are already constructed, we define the root vectors $e_\gamma$ and $f_\gamma$ as
\begin{equation*}
\bm e_\gamma \propto \ldbr \bm e_\alpha, \, \bm e_\beta \rdbr = \bm e_\alpha \bm e_\beta - (-1)^{[\alpha][\beta]} q^{-\bil{\alpha}{\beta}} \bm e_\beta \bm e_\alpha, \qquad \bm f_\gamma \propto \ldbr \bm f_\beta, \, \bm f_\alpha \rdbr = \bm f_\beta \bm f_\alpha - (-1)^{[\alpha][\beta]} q^{\bil{\alpha}{\beta}} \bm f_\alpha \bm f_\beta.
\end{equation*}
The proportionality coefficients are chosen for convenience.

We use the normal order of the positive roots of $\hlslto$ defined as follows. First, we assume that $\alpha_{i j} \prec \alpha_{k l}$ if $i < k$, or if $i = k$ and $j < l$, and $\delta - \alpha_{i j} \prec \delta - \alpha_{k l}$ if $i > k$, or if $i = k$ and $j < l$. Further, we assume that $\alpha_{i j} + m \delta \prec \alpha_{k l} + n \delta$ if $\alpha_{i j} \prec \alpha_{k l}$, or if $\alpha_{i j} = \alpha_{k l}$ and $m < n$, and $(\delta - \alpha_{i j}) + m \delta \prec (\delta - \alpha_{k l}) + n \delta$ if $(\delta - \alpha_{i j}) \prec (\delta - \alpha_{k l})$, or if $(\delta - \alpha_{i j}) = (\delta - \alpha_{k l})$ and $m > n$. Recall that we also assume that the relation (\ref{akd}) is valid. 

The root vectors are defined inductively. We start with the root vectors corresponding to the roots $\alpha_i$ and $-\alpha_i$, which we identify with the generators $\bm e_i$ and $\bm f_i$  of $\uqlslto$:
\begin{gather*}
\bm e_{\alpha_1} = \bm e_{\alpha_{1 2}} = \bm e_1, \qquad \bm e_{\alpha_2} = \bm e_{\alpha_{2 3}} = \bm e_2, \qquad \bm f_{\alpha_1} = \bm f_{\alpha_{1 2}} = \bm f_1, \qquad \bm f_{\alpha_2} = \bm f_{\alpha_{2 3}} = \bm f_2. 
\end{gather*}
It is natural to define \cite{Jim86a}
\begin{equation*}
\bm e_{\alpha_{1 3}} = \ldbr \bm e_{\alpha_{1 2}}, \, \bm e_{\alpha_{2 3}}\rdbr, \qquad \bm f_{\alpha_{1 3}} = \ldbr \bm f_{\alpha_{2 3}}, \, \bm f_{\alpha_{1 2}} \rdbr.
\end{equation*}
Taking into account that $\alpha_0 = \delta - \alpha_{1 3}$, we set
\begin{equation*}
\bm e_{\delta - \alpha_{1 3}} = \bm e_0, \qquad \bm f_{\delta - \alpha_{1 3}} = \bm f_0,
\end{equation*}
and define
\begin{align*}
& \bm e_{\delta - \alpha_{1 2}} = \ldbr \bm e_{\alpha_{2 3}}, \, \bm e_{\delta - \alpha_{1 3}} \rdbr, && \bm e_{\delta - \alpha_{2 3}} = \ldbr \bm e_{\alpha_{1 2}}, \, \bm e_{\delta - \alpha_{1 3}} \rdbr, \\
& \bm f_{\delta - \alpha_{1 2}} = \ldbr \bm f_{\delta - \alpha_{1 3}}, \, \bm f_{\alpha_{2 3}} \rdbr, && \bm f_{\delta - \alpha_{2 3}} = \ldbr \bm f_{\delta - \alpha_{1 3}}, \,  \bm f_{\alpha_{1 2}} \rdbr.
\end{align*}
The root vectors corresponding to the roots $\delta$ and $-\delta$ are indexed by the simple roots $\alpha_i$ and are defined by the relations
\begin{align}
& \bm e'_{\delta, \, \alpha_1} = \ldbr \bm e_{\alpha_1}, \, \bm e_{\delta - \alpha_1} \rdbr, &&
\bm e'_{\delta, \, \alpha_2} = - \ldbr \bm e_{\alpha_2}, \, \bm e_{\delta - \alpha_2} \rdbr, \label{epda} \\
& \bm f'_{\delta, \, \alpha_1} = \ldbr \bm f_{\delta - \alpha_1}, \, \bm f_{\alpha_1} \rdbr, &&
\bm f'_{\delta, \, \alpha_2} = - \ldbr \bm f_{\delta - \alpha_2}, \, \bm f_{\alpha_2} \rdbr.
\end{align}
The remaining definitions for the root vectors corresponding to the real roots are as follows:\footnote{There is no strict rule for determining the root vectors associated with roots containing compound roots of $\Delta_+$. We construct them in a way that allows us to achieve the desired final result.}
\begin{align}
& \bm e_{\alpha_{1 2} + n \delta} = [2]_q^{-1} \ldbr \bm e_{\alpha_{1 2} + (n - 1) \delta}, \, \bm e'_{\delta, \, \alpha_1} \rdbr, \label{ea12pnd} \\
& \bm e_{\alpha_{1 3} + n \delta} = \ldbr \bm e_{\alpha_{1 3} + (n - 1) \delta}, \, \bm e'_{\delta, \, \alpha_2} \rdbr, \label{ea13pnd} \\
& \bm e_{\alpha_{2 3} + n \delta} = - \ldbr \bm e_{\alpha_{2 3} + (n - 1) \delta}, \, \bm e'_{\delta, \, \alpha_1} \rdbr, \label{ea23pnd} \\
& \bm f_{\alpha_{1 2} + n \delta} = [2]_q^{-1} \ldbr \bm f'_{\delta, \, \alpha_1}, \, \bm f_{\alpha_{1 2} + (n - 1) \delta} \rdbr, \\
& \bm f_{\alpha_{1 3} + n \delta} = \ldbr \bm f'_{\delta, \, \alpha_2}, \, \bm f_{\alpha_{1 3} + (n - 1) \delta} \rdbr, \\
& \bm f_{\alpha_{2 3} + n \delta} = - \ldbr \bm f'_{\delta, \, \alpha_1}, \, \bm f_{\alpha_{2 3} + (n - 1) \delta} \rdbr, \\
& \bm e_{(\delta - \alpha_{1 2}) + n \delta} = [2]_q^{-1} \ldbr \bm e'_{\delta, \, \alpha_1}, \, \bm (e_{(\delta - \alpha)_{1 2} + (n - 1) \delta} \rdbr, \label{edma12pnd} \\
& \bm e_{(\delta - \alpha_{1 3}) + n \delta} = \ldbr \bm e'_{\delta, \, \alpha_2}, \, \bm e_{(\delta - \alpha_{1 3}) + (n - 1) \delta} \rdbr, \label{edma13pnd} \\
& \bm e_{(\delta - \alpha_{2 3}) + n \delta} = - \ldbr \bm e'_{\delta, \, \alpha_1}, \, \bm e_{(\delta - \alpha_{2 3}) + (n - 1) \delta} \rdbr, \label{edma23pnd} \\
& \bm f_{(\delta - \alpha_{1 2}) + n \delta} = [2]_q^{-1} \ldbr\bm f_{(\delta - \alpha_{1 2}) + (n - 1) \delta}, \,  \bm f'_{\delta, \, \alpha_1} \rdbr, \\
& \bm f_{(\delta - \alpha_{1 3}) + n \delta} =  \ldbr \bm f_{(\delta - \alpha_{1 3}) + (n - 1) \delta}, \, \bm f'_{\delta, \, \alpha_2} \rdbr, \\
& \bm f_{(\delta - \alpha_{2 3}) + n \delta} = - \ldbr \bm f_{(\delta - \alpha_{2 3}) + (n - 1) \delta}, \,  \bm f'_{\delta, \, \alpha_1} \rdbr,
\end{align}
and corresponding to the imaginary roots are
\begin{align}
& \bm e'_{n \delta, \, \alpha_1} = \ldbr\bm e_{\alpha_1 + (n - 1)\delta}, \, \bm e_{\delta - \alpha_1} \rdbr && \bm e'_{n \delta, \, \alpha_2} = \ldbr\bm e_{\alpha_2 + (n - 1)\delta}, \, \bm e_{\delta - \alpha_2} \rdbr, \label{epnda} \\
& \bm f'_{n \delta, \, \alpha_1} = \ldbr \bm f_{\delta - \alpha_1}, \,  \bm f_{\alpha_1 + (n - 1)\delta} \rdbr && \bm f'_{n \delta, \, \alpha_2} = \ldbr \bm f_{\delta - \alpha_2} , \, \bm f_{\alpha_2 + (n - 1)\delta} \rdbr. \label{fpnda}
\end{align}
The prime in the notation for the root vectors corresponding to the imaginary roots $n \delta$ and $- n \delta$, $n \in \bbZ_{> 0}$, is explained by the fact that to construct the expression for the universal $R$-matrix one uses another set of root vectors corresponding to these roots. They are introduced by the functional equations
\begin{gather}
- \kappa_q \, \bm e_{\delta, \, \alpha_i}(\zeta) = \log(1 - \kappa_q \, \bm e'_{\delta, \, \alpha_i}(\zeta)), \qquad \kappa_q \, \bm f_{\delta, \, \alpha_i}(\zeta^{-1}) = \log(1 + \kappa_q \, \bm f'_{\delta, \, \alpha_i}(\zeta^{-1})), \label{edafda}
\end{gather}
where
\begin{align*}
& \bm e'_{\delta, \, \alpha_i}(\zeta) = \sum_{n = 1}^\infty \bm e'_{n \delta, \, \alpha_i} \, \zeta^n, && \bm e_{\delta, \, \alpha_i}(\zeta) 
= \sum_{n = 1}^\infty \bm e_{n \delta, \, \alpha_i} \, \zeta^n, \\
& \bm f'_{\delta, \, \alpha_i}(\zeta^{-1}) = \sum_{n = 1}^\infty \bm f'_{n \delta, \, \alpha_i} \, \zeta^{- n}, && \bm f_{\delta, \, \alpha_i}(\zeta^{-1}) = \sum_{n = 1}^\infty \bm f_{n \delta, \, \alpha_i} \, \zeta^{- n}.
\end{align*}

\subsection{\texorpdfstring{Universal $R$-matrix}{Universal R-matrix}} \label{s:3.5}

Let $\sigma$ be an automorphism of $\uqlslto \otimes \uqslto$ defined by the equation
\begin{equation*}
\sigma(a_1 \otimes a_2) = (-1)^{[a_1] [a_2]} a_2 \otimes a_1,
\end{equation*}
for all homogeneous $a_1, a_2$. As with any Hopf superalgebra, the quantum superalgebra $\uqlslto$ possesses another comultiplication $\Delta'$, called the opposite comultiplication, defined by the equation
\begin{equation*}
\Delta' = \Delta \circ \sigma.
\end{equation*}
The quantum superalgebra $\uqlslto$ is quasitriangular. It means that there exists a unique invertible even element $\calR \in \uqlslto \otimes \uqlslto$, called the universal $R$-matrix, such that
\begin{equation*}
\Delta'(a) = \calR \, \Delta(a) \, \calR^{-1}
\end{equation*}
for all $a \in \uqlslto$, and
\begin{equation*}
(\Delta \otimes \id) (\calR) = \calR^{(1 3)} \calR^{(2 3)}, \qquad (\id \otimes \Delta) (\calR) = \calR^{(1 3)} \calR^{(1 2)}.
\end{equation*}
Here, if $\calR = \sum_\alpha a_\alpha \otimes b_\alpha$, then
\begin{equation*}
\calR^{(1 2)} = \sum_\alpha a_\alpha \otimes b_\alpha \otimes 1, \qquad \calR^{(1 3)} = \sum_\alpha a_\alpha \otimes 1 \otimes b_\alpha, \qquad \calR^{(2 3)} = \sum_\alpha 1 \otimes a_\alpha \otimes b_\alpha.
\end{equation*}
The most significant property of the universal $R$-matrix is the equality
\begin{equation}
\calR^{(1 2)} \calR^{(1 3)} \calR^{(2 3)} = \calR^{(2 3)} \calR^{(1 3)} \calR^{(1 2)}, \label{rrr}
\end{equation}
known as the Yang--Baxter equation for the universal $R$-matrix.

The expression for the universal $R$-matrix of $\uqlslto$ can be constructed using the procedure proposed by Khoroshkin and Tolstoy~\cite{KhoTol91, KhoTol93a, KhoTol94a, Raz23}. It has the form
\begin{equation*}
\calR = \calR_{\prec \delta} \, \calR_{\sim \delta} \, \calR_{\succ \delta} \, \calK,
\end{equation*}
where $\calR_{\prec \delta}$, $\calR_{\sim \delta}$, $\calR_{\succ \delta}$ and $\calK$ are elements of $\uqlslto \otimes \uqlslto$ described as follows.

The elements $\calR_{\prec \delta}$ and $\calR_{\succ \delta}$ are the products over the sets of roots $\gamma = \alpha_{i j} + n \delta$ and $\gamma = (\delta - \alpha_{i j}) + n \delta$, respectively, of the $q$-exponentials
\begin{equation*}
\calR_\gamma = \exp_{q_\gamma} \big( - \kappa_q (-1)^{[\gamma]} \, a_\gamma^{-1} (\bm e_\gamma^{} \otimes \bm f_\gamma^{}) \big),
\end{equation*}
where
\begin{equation*}
q_\gamma = (-1)^{[\gamma]} q^{(\gamma | \gamma)}
\end{equation*}
and the quantity $a_\gamma$ is determined by the equation
\begin{equation*}
\ldbr \bm e_\gamma, \, \bm f_\gamma \rdbr = a_\gamma \frac{q^{\bm h_\gamma} - q^{- \bm h_\gamma}}{q - q^{-1}}.
\end{equation*}
Here for $\gamma = \sum_{i = 0}^2 m_i \alpha_i$ we define
\begin{equation*}
\bm h_\gamma = \sum_{i = 0}^2 (-1)^{[i]} m_i \bm h_i.
\end{equation*}
The order of the factors in $\calR_{\prec \delta}$ and $\calR_{\succ \delta}$ corresponds to the chosen normal ordering of the roots $\alpha_{i j} + n \delta$ and $(\delta - \alpha_{i j}) + n \delta$.

It is useful to note the following relations:
\begin{gather*}
(-1)^{[\alpha_{i j} + n \delta]} = (-1)^{[(\delta - \alpha_{i j}) + n \delta]} = (-1)^{[i] + [j]}, \\ \bil {(\alpha_{i j}) + n \delta}{(\alpha_{i j}) + n \delta} = \bil {(\delta - \alpha_{i j}) + n \delta}{(\delta - \alpha_{i j}) + n \delta} = (-1)^{[i]} + (-1)^{[j]},
\end{gather*}
which lead to
\begin{gather*}
\\
\bil{\alpha_{1 3}}{\alpha_{1 3}} = \bil{\delta - \alpha_{1 3}}{\delta - \alpha_{1 3}} = 0, \qquad \bil{\alpha_{2 3}}{\alpha_{2 3}} = \bil{\delta - \alpha_{2 3}}{\delta - \alpha_{2 3}} = 0, \\
q_{\alpha_{1 2} + n \delta} = q_{(\delta - \alpha_{1 2}) + n \delta} = q^2, \qquad q_{\alpha_{1 3} + n \delta} = q_{(\delta - \alpha_{1 3}) + n \delta} = -1, \qquad q_{\alpha_{2 3} + n \delta} = q_{(\delta - \alpha_{2 3}) + n \delta} = -1, \\ 
\bil{\alpha_{1 2}}{\alpha_{1 2}} = \bil{\delta - \alpha_{1 2}}{\delta - \alpha_{1 2}} = 2.
\end{gather*}
For the quantities $a_\gamma$, we have
\begin{equation*}
a_{\alpha_{i j} + n \delta} = (-1)^n, \qquad a_{(\delta - \alpha_{i j}) + n \delta} = -(-1)^n.
\end{equation*}
Using these relations, we obtain the following expressions for $\calR_{\prec \delta}$ and $\calR_{\succ \delta}$:
\begin{align*}
& \calR_{\prec \delta} = \overrightarrow{\prod_n} \, \calR_{\alpha_{1 2} + n \delta} \, \overrightarrow{\prod_n} \, \calR_{\alpha_{1 3} + n \delta} \, \overrightarrow{\prod_n} \, \calR_{\alpha_{2 3} + n \delta} \\*
& \hspace{10em} = \overrightarrow{\prod_n} \exp_{q^2}(-\kappa_q(-1)^n (\bm e_{\alpha_{12} + n \delta} \otimes \bm f_{\alpha_{12} + n \delta})) \\*
& \hspace{15em} \times \overrightarrow{\prod_n} \exp_{-1}(\kappa_q(-1)^n (\bm e_{\alpha_{1 3} + n \delta} \otimes \bm f_{\alpha_{1 3} + n \delta})) \\*
& \hspace{20em} \times \overrightarrow{\prod_n} \exp_{-1}(\kappa_q(-1)^n \bm e_{\alpha_{2 3} + n \delta} \otimes \bm f_{\alpha_{2 3} + n \delta})), \\
& \calR_{\succ \delta} = \overrightarrow{\prod_n} \, \calR_{(\delta - \alpha_{2 3}) + n \delta} \, \overrightarrow{\prod_n} \, \calR_{(\delta - \alpha_{1 2}) + n \delta} \, \overrightarrow{\prod_n} \, \calR_{(\delta - \alpha_{1 3}) + n \delta} \\
& \hspace{7em} = \overrightarrow{\prod_n} \exp_{-1}(- \kappa_q (-1)^n (\bm e_{(\delta - \alpha_{2 3}) + n \delta} \otimes \bm f_{(\delta - \alpha_{2 3}) + n \delta})) \\
& \hspace{12em} \times \overrightarrow{\prod_n} \exp_{q^2}(\kappa_q(-1)^n (\bm e_{(\delta - \alpha_{1 2}) + n \delta} \otimes \bm f_{(\delta - \alpha_{1 2}) + n \delta})) \\
& \hspace{17em} \times \overrightarrow{\prod_n} \exp_{-1}(- \kappa_q(-1)^n (\bm e_{(\delta - \alpha_{1 3}) + n \delta} \otimes \bm f_{(\delta - \alpha_{1 3}) + n \delta})).
\end{align*}
The element $\calR_{\sim \delta}$ is defined as
\begin{equation}
\calR_{\sim \delta} = \exp \big( - \kappa_q \sum_{n \in \bbZ_{>0}} \, \sum_{i, j = 1}^2 (-1)^n (u_n)_{i j} \, (\bm e_{n \delta, \, \alpha_i} \otimes \bm f_{n \delta, \, \alpha_j}) \big), \label{rsd}
\end{equation}
where for each $n \in \bbZ_{> 0}$ the quantities $(u_n)_{i j}$ are the matrix entries of the matrix
\begin{equation*}
u_n = \frac{n}{[n]_q} \left( \begin{array}{cc}
0  & - 1 \\
-1 & - [2]_{q^n}
\end{array} \right).
\end{equation*}
Finally, the factor $\calK$ is expressed as
\begin{equation}
\calK = \exp \big(-\hbar \sum_{i, j = 1}^2 c_{i j} (\bm h_i \otimes \bm h_j) \big),
\label{k}
\end{equation}
where $c_{ij}$ are the matrix entries of the matrix
\begin{equation*}
C = A^{-1} = \left( \begin{array}{cc}
0 & -1 \\
-1 & -2
\end{array} \right).
\end{equation*}

\section{Basic monodromy operator}

\subsection{Integrability objects} \label{s:4.1}

To address integrability objects, we introduce the concept of a spectral parameter. This is achieved by endowing $\uqlslto$ with a $\bbZ$-grading. The procedure follows the method described in the paper \cite{BooGoeKluNirRaz14a}. For a given $\zeta \in \bbC^\times$, we define an automorphism $\Gamma_\zeta$ of $\uqlslto$ by its action on the generators of $\uqlslto$ as
\begin{equation*}
\Gamma_\zeta(\bm h_i) = 0, \qquad \Gamma_\zeta(\bm e_i) = \zeta^{s_i} \bm e_i, \qquad \Gamma_\zeta(\bm f_i) = \zeta^{-s_i} \bm f_i,
\end{equation*}
where $s_i$ are arbitrary integers. The automorphisms $\Gamma_\zeta$ induce the $\bbZ$-grading on $\uqlslto$ with the grading subspaces
\begin{equation*}
\uqlslto_m = \{ a \in \uqlslto \mid \Gamma_\zeta(a) = \zeta^m a \}.
\end{equation*}
It is noteworthy that for any $\zeta \in \bbC^\times$ the universal $R$-matrix of $\uqlslto$ satisfies the condition
\begin{equation}
(\Gamma_\zeta \otimes \Gamma_\zeta)(\calR) = \calR. \label{phiphir}
\end{equation}
We use the following notations:
\begin{equation*}
s = s_0 + s_1 + s_2
\end{equation*}
and
\begin{equation*}
s_{1 2} = s_1, \qquad s_{1 3} = s_1 + s_2, \qquad s_{2 3} = s_2.
\end{equation*}
In what follows, we assume that $s \ne 0$.

For any homomorphism $\varphi$ from $\uqlslto$ to an associative superalgebra $A$, we define a family of homomorphisms $\varphi_\zeta$, $\zeta \in \bbC^\times$, from $\uqlslto$ to $A$ by the relation
\begin{equation*}
\varphi_\zeta = \varphi \circ \Gamma_\zeta.
\end{equation*}
Let $\varphi_1$ and $\varphi_2$ be homomorphisms from $\uqlslto$ to associative superalgebras $A_1$ and $A_2$. The corresponding integrability object $X_{\varphi_1 | \varphi_2}(\zeta_1 | \zeta_2)$ is defined by the equation
\begin{equation*}
X_{\varphi_1 | \varphi_2}(\zeta_1 | \zeta_2) = \rho_{\varphi_1 | \varphi_2}(\zeta_1 | \zeta_2) (\varphi_{\zeta_1} \otimes \varphi_{\zeta_2})(\calR),
\end{equation*}
where $\zeta_1$ and $\zeta_2$ are spectral parameters, and $\rho_{\varphi_1 | \varphi_2}(\zeta_1 | \zeta_2)$ is a scalar normalization factor. By definition, $X_{\varphi_1 | \varphi_2}(\zeta_1 | \zeta_2)$ is an element of the superalgebra $A_1 \otimes A_2$.

\subsection{Monodromy operators} \label{s:4.2}

In the present paper, we study a special class of integrability objects known as monodromy operators. To construct these operators, we use either an infinite-dimensional representation $\widetilde{\varphi}^\lambda$ or a finite-dimensional representation  $\varphi^\lambda$ for the first factor of $\calR$. These representations were defined in section \ref{s:3.3}. For the second factor, we use the finite-dimensional representation $\varphi = \varphi^{(1,0,0)}$. Thus, we consider the monodromy operators
\begin{equation*}
\widetilde M^\lambda(\zeta | \eta) = ((\widetilde \varphi^\lambda)_\zeta \otimes \varphi_\eta)(\calR), \qquad M^\lambda(\zeta | \eta) = ((\varphi^\lambda)_\zeta \otimes \varphi_\eta)(\calR).
\end{equation*}
In fact, it is convenient for applications to define the basic monodromy operator
\begin{equation*}
M(\zeta | \eta) = (\epsilon_\zeta \otimes \varphi_\eta))(\calR).
\end{equation*}
For any $\zeta$ and $\eta$ the operator $M(\zeta | \eta)$ is an element of the space $\uqglto \otimes \End(V)$. The monodromy operators $\widetilde M^\lambda(\zeta | \eta)$ and $M^\lambda(\zeta | \eta)$ can be constructed with the help of the relations
\begin{equation*}
\widetilde M^\lambda(\zeta | \eta) = (\widetilde \pi^\lambda \otimes \id)(M(\zeta | \eta)), \qquad M^\lambda(\zeta | \eta) = (\pi^\lambda \otimes \id)(M(\zeta | \eta)).
\end{equation*}

Using equation (\ref{phiphir}), we demonstrate that
\begin{equation}
M(\zeta \nu | \eta \nu) = M(\zeta | \eta) \label{mznen}
\end{equation}
for any $\nu \in \bbC^\times$. We have
\begin{equation*}
M(\zeta | \eta) = M(\zeta \eta^{-1}),
\end{equation*}
where $M(\zeta) = M(\zeta | 1)$.

The monodromy operator $M^\lambda(\zeta_1 | \zeta_2)$ for $\lambda = (1, \, 0, \, 0)$ plays a special role and is referred to as the $R$-operator. We denote it as $R(\zeta_1 | \zeta_2)$. It follows from equation (\ref{rrr}) that $R(\zeta_1 | \zeta_2)$ satisfies the Yang--Baxter equation:
\begin{equation*}
R^{(1 2)}(\zeta_1 | \zeta_2) R^{(1 3)}(\zeta_1 | \zeta_3) R^{(2 3)}(\zeta_2 | \zeta_3) = R^{(2 3)}(\zeta_2 | \zeta_3) R^{(1 3)}(\zeta_1 | \zeta_3) R^{(1 2)}(\zeta_1 | \zeta_2).
\end{equation*}
Here, if $\calR = \sum_\alpha a_\alpha \otimes b_\alpha$, then
\begin{gather*}
R^{(1 2)}(\zeta_1 | \zeta_2)  = \sum_\alpha \varphi_{\zeta_1}(a_\alpha) \otimes \varphi_{\zeta_2}(b_\alpha) \otimes \id_V, \quad R^{(1 3)}(\zeta_1| \zeta_3) = \sum_\alpha \varphi_{\zeta_1}(a_\alpha) \otimes \id_V \otimes \varphi_{\zeta_3}(b_\alpha), \\
R^{(2 3)}(\zeta_2 | \zeta_3) = \sum_\alpha \id_V \otimes \varphi_{\zeta_2}(a_\alpha) \otimes \varphi_{\zeta_3}(b_\alpha).
\end{gather*}
Applying the mapping $\epsilon_\zeta \otimes \varphi_{\zeta_1} \otimes \varphi_{\zeta_2}$, to both sides of equation (\ref{rrr}), we obtain
\begin{equation}
M^{(1 2)}(\zeta | \zeta_1) M^{(1 3)}(\zeta | \zeta_2) R^{(2 3)}(\zeta_1 | \zeta_2) = R^{(2 3)}(\zeta_1 | \zeta_2) M^{(1 3)}(\zeta | \zeta_2) M^{(1 2)}(\zeta | \zeta_1). \label{mmr}
\end{equation}
Here, if $\calR = \sum_\alpha a_\alpha \otimes b_\alpha$, then
\begin{gather*}
M^{(1 2)}(\zeta | \zeta_1)  = \sum_\alpha \epsilon_{\zeta}(a_\alpha) \otimes \varphi_{\zeta_1}(b_\alpha) \otimes \id_V, \quad M^{(1 3)}(\zeta| \zeta_2) = \sum_\alpha \epsilon_{\zeta}(a_\alpha) \otimes \id_V \otimes \varphi_{\zeta_2}(b_\alpha), \\
R^{(2 3)}(\zeta_1 | \zeta_2) = \sum_\alpha \id_V \otimes \varphi_{\zeta_1}(a_\alpha) \otimes \varphi_{\zeta_2}(b_\alpha).
\end{gather*}
Substituting $\zeta = \zeta_1 \zeta_2$ in (\ref{mmr}) and using equation (\ref{mznen}), we arrive at the more familiar equation
\begin{equation*}
 R^{(2 3)}(\zeta_1 | \zeta_2) M^{(1 3)}(\zeta_1) M^{(1 2)}(\zeta_2) = M^{(1 2)}(\zeta_2) M^{(1 3)}(\zeta_1) R^{(2 3)}(\zeta_1 | \zeta_2). \label{rmm}
\end{equation*}

\subsection{Structure of basic monodromy operator} \label{s:4.3}

Let us first discuss the general structure of the monodromy operator $M(\zeta)$. We have
\begin{equation}
M(\zeta) = U(\zeta) V(\zeta) W(\zeta) D(\zeta), \label{muvwd}
\end{equation}
where
\begin{gather*}
U(\zeta) = (\epsilon_\zeta \otimes \varphi)(\calR_{\prec \delta}), \qquad V(\zeta) = (\epsilon_\zeta \otimes \varphi)(\calR_{\sim \delta}), \qquad
W(\zeta) = (\epsilon_\zeta \otimes \varphi)(\calR_{\succ \delta}), \\
D(\zeta) = (\epsilon_\zeta \otimes \varphi)(\calK).
\end{gather*}

It is not difficult to find the explicit form of $D(\zeta)$. Note that although it does not actually depend on $\zeta$, we write $D(\zeta)$ for consistency. Using equations (\ref{k}), (\ref{j2h}), and (\ref{fh}), we obtain
\begin{equation}
D(\zeta) = q^{\bm K_2 + \bm K_3} \otimes \bbM_{1 1} + q^{\bm K_1 + \bm K_3} \otimes \bbM_{2 2} + q^{\bm K_1 + \bm K_2 + 2 \bm K_3} \otimes \bbM_{3 3}. \label{dz}
\end{equation}

It is clear that
\begin{equation*}
\epsilon_\zeta(\bm e_{\alpha_{i j} + n \delta}) = \zeta^{s_{i j} + n s} \epsilon(\bm e_{\alpha_{i j} + n \delta}).
\end{equation*}
Using this relation and the equalities \cite{Raz23}
\begin{align*}
& \varphi(\bm f_{\alpha_{1 2} + n \delta}) = (-1)^n q^{-2 n} \bbM_{2 1}, && \varphi(\bm f_{(\delta - \alpha_{ 1 2}) + n \delta}) = (-1)^n q^{- 2n  - 1} \bbM_{1 2}, \\
& \varphi(\bm f_{\alpha_{1 3} + n \delta}) = q^{- 2 n} \bbM_{3 1}, && \varphi(\bm f_{(\delta - \alpha_{1 3}) + n \delta}) = q^{- 2 n - 1} \bbM_{1 3}, \\
& \varphi(\bm f_{\alpha_{2 3} + n \delta}) = (-1)^n q^{- 3n} \bbM_{3 2}, && \varphi(\bm f_{(\delta - \alpha_{2 3}) + n \delta}) = - (-1)^n q^{- 3n - 2} \bbM_{2 3},
\end{align*}
we get
\begin{align*}
(\epsilon_\zeta \otimes \varphi) \big( \, \overrightarrow{\prod_n} \, \calR_{\alpha_{1 2} + n \delta} \big) & = \overrightarrow{\prod_n} \exp_{q^2}(-\kappa_q(-1)^n (\epsilon_\zeta(\bm e_{\alpha_{12} + n \delta}) \otimes \varphi(\bm f_{\alpha_{12} + n \delta}))) \\
& = \overrightarrow{\prod_n} \exp_{q^2}(-\kappa_q \zeta^{s_{1 2}} ((q^{-2} \zeta^s)^n \epsilon(\bm e_{\alpha_{12} + n \delta}) \otimes \bbM_{2 1})) \\
& = \overrightarrow{\prod_n} (1 \otimes \mathbbm 1 - \kappa_q \zeta^{s_{1 2}} (q^{-2} \zeta^s)^n \epsilon(\bm e_{\alpha_{12} + n \delta}) \otimes \bbM_{2 1}) \\
& = 1 \otimes \mathbbm 1 - \kappa_q \zeta^{s_{1 2}} \sum_n (q^{-2} \zeta^s)^n \epsilon(\bm e_{\alpha_{12} + n \delta}) \otimes \bbM_{2 1} \\
& = 1 \otimes \mathbbm 1 - \kappa_q \zeta^{s_{1 2}}  \bm{\mathcal E}_{\alpha_{1 2}}(q^{-2} \zeta^s) \otimes \bbM_{2 1}.
\end{align*}
Here and below, the notation $\mathbbm 1$ is used for the unit element of $\End(V)$ and, moreover, we define the following generating functions: 
\begin{equation*}
\bm{\mathcal E}_{\alpha_{i j}}(\zeta) = \epsilon(\bm e_{\alpha_{i j}}(\zeta)) = \sum_{n = 0}^\infty \epsilon(\bm e_{\alpha_{i j} + n \delta}) \zeta^n, \qquad \bm{\mathcal E}_{\delta - \alpha_{i j}}(\zeta) = \epsilon(\bm e_{\delta - \alpha_{i j}}(\zeta)) = \sum_{n = 0}^\infty \epsilon(\bm e_{(\delta - \alpha_{i j}) + n \delta}) \zeta^n.
\end{equation*}
Similarly, we obtain
\begin{align*}
& (\epsilon_\zeta \otimes \varphi) \big( \, \overrightarrow{\prod_n} \, \calR_{\alpha_{1 3} + n \delta} \big) = 1 \otimes \mathbbm 1 + \kappa_q \zeta^{s_{1 3}} \bm{\mathcal E}_{\alpha_{1 3}}(- q^{-2} \zeta^s) \otimes \bbM_{3 1}, \\
& (\epsilon_\zeta \otimes \varphi) \big( \, \overrightarrow{\prod_n} \, \calR_{\alpha_{2 3} + n \delta} \big) = 1 \otimes \mathbbm 1 + \kappa_q \zeta^{s_{2 3}} \bm{\mathcal E}_{\alpha_{2 3}}(q^{-3} \zeta^s) \otimes \bbM_{3 2}, \\
& (\epsilon_\zeta \otimes \varphi) \big( \, \overrightarrow{\prod_n} \, \calR_{(\delta - \alpha_{1 2}) + n \delta} \big) = 1 \otimes \mathbbm 1 + \kappa_q \zeta^{s - s_{1 2}} q^{-1} \bm{\mathcal E}_{\delta - \alpha_{1 2}}(q^{-2} \zeta^s) \otimes \bbM_{1 2}, \\
& (\epsilon_\zeta \otimes \varphi) \big( \, \overrightarrow{\prod_n} \, \calR_{(\delta - \alpha_{1 3}) + n \delta} \big) = 1 \otimes \mathbbm 1 + \kappa_q \zeta^{s - s_{1 3}} q^{-1} \bm{\mathcal E}_{\delta - \alpha_{1 3}}(-q^{-2} \zeta^s) \otimes \bbM_{1 3}, \\
& (\epsilon_\zeta \otimes \varphi) \big( \, \overrightarrow{\prod_n} \, \calR_{(\delta - \alpha_{2 3}) + n \delta} \big) = 1 \otimes \mathbbm 1 - \kappa_q \zeta^{s - s_{2 3}} q^{-2} \bm{\mathcal E}_{\delta - \alpha_{2 3}}(q^{-3} \zeta^s) \otimes \bbM_{2 3}.
\end{align*}
Representing $U(\zeta)$ and $W(\zeta)$ as
\begin{equation*}
U(\zeta) = 1 \otimes \mathbbm 1 + \sum_{\substack{i, j = 1 \\ i > j}}^3 U_{i j}(\zeta) \otimes \bbM_{i j}, \qquad W(\zeta) = 1 \otimes \mathbbm 1 + \sum_{\substack{i, j = 1 \\ i < j}}^3 W_{i j}(\zeta) \otimes \bbM_{i j},
\end{equation*}
we see that
\begin{align}
& U_{2 1}(\zeta) = - \kappa_q \zeta^{s_{1 2}} \bm \calE_{\alpha_{1 2}}(q^{-2} \zeta^s), \label{u21} \\
& U_{3 1}(\zeta) = \kappa_q \zeta^{s_{1 3}} \bm \calE_{\alpha_{1 3}}(-q^{-2} \zeta^s), \label{u31} \\
& U_{3 2}(\zeta) = \kappa_q \zeta^{s_{2 3}} \bm \calE_{\alpha_{2 3}}(q^{-3} \zeta^s), \label{u32} \\
& W_{1 2}(\zeta) = \kappa_q \zeta^{s - s_{1 2}} q^{-1} \bm \calE_{\delta - \alpha_{1 2}}(q^{-2} \zeta^s), \\
& W_{1 3}(\zeta) = \kappa_q \zeta^{s - s_{1 3}} q^{-1} \bm \calE_{\delta - \alpha_{1 3}}(-q^{-2} \zeta^s), \\
& W_{2 3}(\zeta) = - \kappa_q \zeta^{s - s_{2 3}} q^{-2} \bm \calE_{\delta - \alpha_{2 3}}(q^{-3} \zeta^s).
\end{align}

It is clear that
\begin{equation*}
\epsilon_\zeta(\bm e_{n \delta, \, \alpha_i}) = \zeta^{n s} \epsilon\bm(e_{n \delta, \, \alpha_i}).
\end{equation*}
Using this relation, equation (\ref{rsd}), and the equalities \cite{Raz23}
\begin{equation*}
\varphi(\bm f_{n \delta, \, \alpha_1}) = (-1)^{- n - 1} q^{-n} \frac{[n]_q}{n} (\bbM_{1 1}  - q^{-2 n} \bbM_{2 2}), \qquad
\varphi(\bm f_{n \delta, \, \alpha_2}) = (-1)^{- n - 1} q^{-2n} \frac{[n]_q}{n} (\bbM_{2 2} + \bbM_{3 3}),
\end{equation*}
we arrive at the expression
\begin{multline*}
V(\zeta) = \exp \big( - \kappa_q \sum_{n = 1}^\infty \zeta^{n s} q^{-2 n} (q^n \epsilon(\bm e_{n \delta, \, \alpha_2}) \otimes \bbM_{1 1} + \\ (\epsilon(\bm e_{n \delta, \, \alpha_1}) + q^n \epsilon(\bm e_{n \delta, \, \alpha_2})) \otimes \bbM_{2 2} + (\epsilon(\bm e_{n \delta, \, \alpha_1}) + [2]_{q^n} \epsilon(\bm e_{n \delta, \, \alpha_2})) \otimes \bbM_{3 3} \big).
\end{multline*}
Representing $V(\zeta)$ as
\begin{equation*}
V(\zeta) = \sum_{i = 1}^3 V_{i i} \otimes \bbM_{i i},
\end{equation*}
we see that
\begin{align}
& V_{1 1}(\zeta) = \exp \big( - \kappa_q \bm \calE_{\delta, \, \alpha_2}(q^{-1} \zeta^s) \big),\label{v11} \\
& V_{2 2}(\zeta) = \exp \big( - \kappa_q \bm \calE_{\delta, \, \alpha_1}(q^{-2} \zeta^s) - \kappa_q \bm \calE_{\delta, \, \alpha_2}(q^{-1} \zeta^s) \big), \label{v22} \\
& V_{3 3}(\zeta) = \exp \big( - \kappa_q \bm \calE_{\delta, \, \alpha_1}(q^{-2} \zeta^s) - \kappa_q \bm \calE_{\delta, \, \alpha_2}(q^{-1} \zeta^s) - \kappa_q \bm \calE_{\delta, \, \alpha_2}(q^{-3} \zeta^s) \big), \label{v33} 
\end{align}
where we use the notation
\begin{equation*}
\bm \calE_{\delta, \, \alpha_i}(\zeta) = \sum_{n = 1}^\infty \epsilon(\bm e_{n \delta, \, \alpha_i}) \zeta^n.
\end{equation*}

\subsection{Ansatz for basic monodromy operator} \label{s:4.4}

We have
\begin{equation*}
R(\zeta) = (\pi \otimes \id)(M(\zeta)).
\end{equation*}
Using this relation, we formulate an ansatz for the basic monodromy operator $M(\zeta)$. The explicit expression for $R(\zeta)$ has the form\footnote{Note that the paper \cite{Raz23} contains a misprint in the expression for the scalar factor of the $R$-operator. The corrected version of the paper \cite{Raz23} can be found at \href{https://arxiv.org/abs/2210.12721v2}{\tt arXiv.org}.}
\begin{multline*}
R(\zeta) = \frac{1}{1 - q^2 \zeta^s} \Big[ (1 - q^2 \zeta^s) \sum_{i = 1}^2  \bbM_{i i} \otimes  \bbM_{i i} + q^2 (1 - q^{-2} \zeta^s)  \bbM_{3 3} \otimes  \bbM_{3 3} + q (1 - \zeta^s) \sum_{\substack{i, j = 1 \\ i \ne j}}^3  \bbM_{i i} \otimes  \bbM_{j j} \\
{} + (1 - q^2) \sum_{\substack{i, j = 1 \\ i < j}}^3 (-1)^{[j]} \, \zeta^{s_{i j}} \,  \bbM_{i j} \otimes  \bbM_{j i} + (1 - q^2) \sum_{\substack{i, j = 1 \\ i > j}}^3 (-1)^{[j]} \, \zeta^{s - s_{j i}} \,  \bbM_{i j} \otimes  \bbM_{j i} \Big].
\end{multline*}
see the paper \cite{Raz23}. We represent $R(\zeta)$ in the factorized form
\begin{equation*}
R(\zeta) = \rme^{f(\zeta^s)} S(\zeta) K(\zeta),
\end{equation*}
where $f(\zeta) = - \log(1 - q^2 \zeta)$, 
\begin{multline}
K(\zeta) = (\varphi_\zeta \otimes \varphi)(\calK) = (\pi \otimes \id) (D(\zeta)) \\=  \bbM_{1 1} \otimes  \bbM_{1 1} +  \bbM_{2 2} \otimes  \bbM_{2 2} + q^2  \bbM_{3 3} \otimes  \bbM_{3 3} + q \sum_{\substack{i, j = 1 \\ i \ne j}}^3  \bbM_{i i} \otimes  \bbM_{j j}, \label{kz}
\end{multline}
and
\begin{multline*}
S(\zeta) = (1 - q^2 \zeta^s)  \bbM_{1 1} \otimes  \bbM_{1 1} + (1 - q^2 \zeta^s)  \bbM_{2 2} \otimes  \bbM_{2 2} + (1 - q^{-2} \zeta^s)  \bbM_{3 3} \otimes  \bbM_{3 3} \\*
+ (1 - \zeta^s) \sum_{\substack{i, j = 1 \\ i \ne j}}^3  \bbM_{i i} \otimes  \bbM_{j j} + \kappa_q \sum_{\substack{i, j = 1 \\ i < j}}^3 (-1)^{[i]} \, \zeta^{s - s_{i j}} \,  \bbM_{j i} \otimes  \bbM_{i j} + \kappa_q \sum_{\substack{i, j = 1 \\ i > j}}^3 (-1)^{[i]} \, \zeta^{s_{j i}} \,  \bbM_{j i} \otimes  \bbM_{i j}.
\end{multline*}
It is convenient to represent $S(\zeta)$ as:
\begin{multline}
S(\zeta) = (\pi \otimes \id) \big( (1 - q^{2 \bm K_1} \zeta^s) \otimes  \bbM_{1 1} + (1 - q^{2 \bm K_2} \zeta^s) \otimes  \bbM_{2 2} + (1 - q^{-2 \bm K_3} \zeta^s) \otimes  \bbM_{3 3} \\*
+ \kappa_q \sum_{\substack{i, j = 1 \\ i < j}}^3 (-1)^{[i]} \, \zeta^{s - s_{i j}} \, \bm F_{i j} \otimes  \bbM_{i j} + \kappa_q \sum_{\substack{i, j = 1 \\ i > j}}^3 (-1)^{[i]} \, \zeta^{s_{j i}} \, \bm E_{j i} \otimes  \bbM_{i j} \big). \label{sz}
\end{multline}
Using (\ref{kz}) and (\ref{sz}), we formulate an ansatz for the basic monodromy operator, which we denote as $N(\zeta)$. We assume that
\begin{equation}
N(\zeta) = O(\zeta) D(\zeta), \label{nod}
\end{equation}
where $D(\zeta)$ is given by equation (\ref{dz}), and for $O(\zeta)$ we use the ansatz
\begin{multline*}
O(\zeta) = (1 - q^{2 \bm K_1} \zeta^s) \otimes  \bbM_{1 1} + (1 - q^{2 \bm K_2} \zeta^s) \otimes  \bbM_{2 2} + (1 - q^{-2 \bm K_3} \zeta^s) \otimes  \bbM_{3 3} \\*
 + \kappa_q \sum_{\substack{i, j = 1 \\ i < j}}^3 c_{i j} \, \zeta^{s - s_{i j}} \, \bm F_{i j} \, q^{\sum_{k = 1}^3 c_{i j k} \bm K_k} \otimes  \bbM_{i j} + \kappa_q \sum_{\substack{i, j = 1 \\ i > j}}^3 d_{i j} \, \zeta^{s_{j i}} \, \bm E_{j i} \, q^{\sum_{k = 1}^3 d_{i j k} \bm K_k} \otimes  \bbM_{i j}.
\end{multline*}
Here $c_{i j}$, $d_{i j}$, $c_{i j k}$ and $d_{i j k}$ are some complex numbers. Substituting (\ref{nod}) into the equation
\begin{equation}
 R^{(2 3)}(\zeta_1 | \zeta_2) N^{(1 3)}(\zeta_1) N^{(1 2)}(\zeta_2) =  N^{(1 2)}(\zeta_2) N^{(1 3)}(\zeta_1) R^{(2 3)}(\zeta_1 | \zeta_2), \label{rnn}
\end{equation}
we see that this equation is satisfied if
\begin{multline*}
O(\zeta) = \big( (1 - q^{2 \bm K_1} \zeta^s) \otimes  \bbM_{1 1} + (1 - q^{2 \bm K_2} \zeta^s) \otimes  \bbM_{2 2} + (1 - q^{-2 \bm K_3} \zeta^s) \otimes  \bbM_{3 3} \\*
- \kappa_q \sum_{\substack{i, j = 1 \\ i < j}}^3 (-1)^{[i]} \, \zeta^{s - s_{i j}} q^{-1} \bm F_{i j} q^{(-1)^{[i]} \bm K_i + \bm K_j} \otimes  \bbM_{i j} - \kappa_q \sum_{\substack{i, j = 1 \\ i > j}}^3 (-1)^{[i]} \, \zeta^{s_{j i}} \, \bm E_{j i} \otimes  \bbM_{i j}.
\end{multline*}
It is clear that solution of the equation (\ref{rnn}) is defined up to a central factor and an automorphism of $\uqglto$. We fix this ambiguity to achieve the desired final result.

Let us demonstrate that for some function $\Phi(\zeta)$ taking values in the center of $\uqglto$, the equation
\begin{equation*}
M(\zeta) = (\rme^{\Phi(\zeta)} \otimes \mathbbm 1) N(\zeta)
\end{equation*}
is satisfied. Taking into account relations (\ref{muvwd}) and (\ref{nod}), we see that the above equation is equivalent to the equation
\begin{equation}
\rme^{\Phi(\zeta)} O(\zeta) = U(\zeta) V(\zeta) W(\zeta). \label{efouvw}
\end{equation}
Rewrite this equation in the component form, and resolve the obtained equations with respect to the components of $U(\zeta)$, $V(\zeta)$ and $W(\zeta)$. We come to the system
\begin{align*}
& \rme^{\Phi(\zeta)} O_{1 1}(\zeta) = V_{1 1}(\zeta), \\
& \rme^{\Phi(\zeta)} O_{1 2}(\zeta) = V_{1 1}(\zeta)(\zeta) W_{1 2}(\zeta), \\
& \rme^{\Phi(\zeta)} O_{1 3}(\zeta) = V_{1 1}(\zeta) W_{1 3}(\zeta), \\
& \rme^{\Phi(\zeta)} O_{2 1}(\zeta) = U_{2 1}(\zeta) V_{1 1}(\zeta), \\
& \rme^{\Phi(\zeta)} O_{2 2}(\zeta) = V_{2 2}(\zeta) + U_{2 1}(\zeta) V_{1 1}(\zeta) W_{1 2}(\zeta), \\
& \rme^{\Phi(\zeta)} O_{2 3}(\zeta) = V_{2 2}(\zeta) W_{2 3}(\zeta) + U_{2 1}(\zeta) V_{1 1}(\zeta) W_{1 3}(\zeta), \\
& \rme^{\Phi(\zeta)} O_{3 1}(\zeta) = U_{3 1}(\zeta) V_{1 1}(\zeta), \\
& \rme^{\Phi(\zeta)} O_{3 2}(\zeta) = U_{3 2}(\zeta) V_{2 2}(\zeta) + U_{3 1}(\zeta) V_{1 1}(\zeta) W_{1 2}(\zeta), \\
& \rme^{\Phi(\zeta)} O_{3 3}(\zeta) = V_{3 3}(\zeta) - U_{3 1}(\zeta) V_{1 1}(\zeta) W_{1 3}(\zeta) - U_{3 2}(\zeta) V_{2 2}(\zeta) W_{2 3}(\zeta) 
\end{align*}
and obtain
\begin{align}
& U_{2 1}(\zeta) = O_{2 1}(\zeta) O_{1 1}(\zeta)^{-1}, \label{u21o} \\
& U_{3 1}(\zeta) = O_{3 1}(\zeta) O_{1 1}(\zeta)^{-1}, \label{u31o} \\
& U_{3 2}(\zeta) = O'_{3 2}(\zeta) O'^{-1}_{2 2}(\zeta), \label{u32o} \\
& V_{1 1}(\zeta) = \rme^{\Phi(\zeta)} O_{1 1}(\zeta), \label{v11o} \\
& V_{2 2}(\zeta) = \rme^{\Phi(\zeta)} O_{2 2}'(\zeta), \label{v22o}\\
& V_{3 3}(\zeta) = \rme^{\Phi(\zeta)} O''_{3 3}(\zeta), \label{v33o} \\
& W_{1 2}(\zeta) = O_{1 1}^{-1}(\zeta) O_{1 2}(\zeta), \label{w12o} \\
& W_{1 3}(\zeta) = O_{1 1}^{-1}(\zeta) O_{1 3}(\zeta), \label{w13o} \\
& W_{2 3}(\zeta) = O'^{-1}_{2 2}(\zeta) O'_{2 3}(\zeta), \label{w23o}
\end{align}
where
\begin{align*}
& O'_{2 2}(\zeta) = O_{2 2}(\zeta) - O_{2 1}(\zeta) O^{-1}_{1 1}(\zeta) O_{1 2}(\zeta), \\
& O'_{2 3}(\zeta) = O_{2 3}(\zeta) - O_{2 1}(\zeta) O_{1 1}^{-1}(\zeta) O_{1 3}(\zeta), \\
& O'_{3 2}(\zeta) = O_{3 2}(\zeta) - O_{3 1}(\zeta) O^{-1}_{1 1}(\zeta) O_{1 2}(\zeta), \\
& O'_{3 3}(\zeta) = O_{33}(\zeta) + O_{3 1}(\zeta) O^{-1}_{1 1}(\zeta) O_{1 3}(\zeta), \\
& O''_{3 3}(\zeta) = O'_{33}(\zeta) + O'_{32}(\zeta) O'^{-1}_{2 2}(\zeta) O'_{2 3}(\zeta).
\end{align*}
It is evident that equations (\ref{u21o})--(\ref{w23o}) are equivalent to equation (\ref{efouvw}).

To prove the validity of equations (\ref{u21o})--(\ref{u32o}), we use the relations
\begin{align}
& \bm \calE_{\alpha_{1 2}}(\zeta) = \epsilon(\bm e_{\alpha_{1 2}}) + \frac{1}{[2]_q} \zeta \, [\bm \calE_{\alpha_{1 2}}(\zeta), \, \epsilon\bm (\bm e'_{\delta, \, \alpha_1})], \label{ea12e} \\
& \bm \calE_{\alpha_{1 3}}(\zeta) = \epsilon(\bm e_{\alpha_{1 3}}) + \zeta \, [\bm \calE_{\alpha_{1 3}}(\zeta), \, \epsilon(\bm e'_{\delta, \, \alpha_2})], \label{ea13e} \\[.5em]
& \bm \calE_{\alpha_{2 3}}(\zeta) = \epsilon(\bm e_{\alpha_{2 3}}) - \zeta \,  [\bm \calE_{\alpha_{2 3}}(\zeta), \, \epsilon(\bm e'_{\delta, \, \alpha_1})], \label{ea23e}
\end{align}
which follow from (\ref{ea12pnd})--(\ref{ea23pnd}). These relations determine the generating functions $\bm \calE_{\alpha_{1 2}}(\zeta)$, $\bm \calE_{\alpha_{1 3}}(\zeta)$ and $\bm \calE_{\alpha_{2 3}}(\zeta)$ uniquely.

We start with the proof of equation (\ref{u21o}). It follows from (\ref{u21}) that this equation is equivalent to
\begin{equation}
\bm \calE_{\alpha_{1 2}}(\zeta) = - \kappa_q^{-1} \zeta^{-s_{1 2}} O_{2 1}(\zeta) O_{1 1}(\zeta)^{-1}|_{\zeta \to q^{2/s} \zeta^{1/s}} \, . \label{ea12oo}
\end{equation}
We get convinced that the right hand side of this equation satisfies equation (\ref{ea12e}).
We see that $\bm \calE_{\alpha_{1 2}}(\zeta)$ and the right hand side of (\ref{ea12oo}) satisfy the same equation. This means that equation (\ref{ea12oo}) is valid and, therefore, equation (\ref{u21o}) is also valid.

We now consider equation (\ref{u31o}). It follows from (\ref{u31}) that this equation is equivalent to
\begin{equation}
\bm \calE_{\alpha_{1 3}}(\zeta) = \kappa_q^{-1} \zeta^{-s_{1 3}} O_{3 1}(\zeta)  O_{1 1}(\zeta)^{-1} |_{\zeta \to (-q^2 \zeta)^{1/s}}. \label{ea13oo}
\end{equation}
We substitute the right-hand side of this equation for $\bm \calE_{\alpha_{1 3}}(\zeta)$ into equation (\ref{ea13e}) and observe that it is satisfied. We note that the right-hand side of equation (\ref{ea13oo}) and $\bm \calE_{\alpha_{1 3}}(\zeta)$ satisfy the same equations, and conclude that equation (\ref{ea13oo}) is true. Hence, equation (\ref{u31o}) is also true.

Finally, we proceed to equation (\ref{u32o}). It follows from equation (\ref{u32}) that it can be written as
\begin{equation}
\bm \calE_{\alpha_{2 3}}(\zeta) = \kappa_q^{-1} \zeta^{-s_{2 3}} O'_{3 2}(\zeta)  O_{2 2}'(\zeta)^{-1} |_{\zeta \to (q^3 \zeta)^{1/s}}. \label{ea23oo}
\end{equation}
We verify that the right hand side of this equation satisfies equation (\ref{ea23e}).
We see that $\bm \calE_{\alpha_{2 3}}(\zeta)$ and the right hand side of (\ref{ea23oo}) satisfy the same equation with the same initial condition. This means that equation (\ref{ea23oo}) is valid and, therefore, equation (\ref{u32o}) is also valid.

In the same way, using the relations
\begin{align*}
& \bm \calE_{\delta - \alpha_{1 2}}(\zeta) = \epsilon(\bm e_{\delta - \alpha_{1 2}}) + \frac{1}{[2]_q} \zeta \, [\bm \calE_{\delta - \alpha_{1 2}}(\zeta), \, \epsilon\bm (\bm e'_{\delta, \, \alpha_1})], \\
& \bm \calE_{\delta - \alpha_{1 3}}(\zeta) = \epsilon(\bm e_{\delta - \alpha_{1 3}}) + \zeta \, [\bm \calE_{\delta - \alpha_{1 3}}(\zeta), \, \epsilon(\bm e'_{\delta, \, \alpha_2})], \\[.5em]
& \bm \calE_{\delta - \alpha_{2 3}}(\zeta) = \epsilon(\bm e_{\delta - \alpha_{2 3}}) + \zeta \,  [\bm \calE_{\delta - \alpha_{2 3}}(\zeta), \, \epsilon(\bm e'_{\delta, \, \alpha_1})],
\end{align*}
which follow from (\ref{edma12pnd})--(\ref{edma23pnd}), we demonstrate that equations (\ref{w12o})--(\ref{w23o}) are true.

It remains to verify equations (\ref{v11o})--(\ref{v33o}). We start with equation (\ref{v11o}). We have
\begin{equation*}
V_{1 1}(\zeta) = \exp(-\kappa \bm \calE_{\delta, \, \alpha_2} (q^{-1} \zeta^s)) = \exp \big( - \kappa_q \sum_{n = 1}^\infty \zeta^{n s} q^{-n} \bm \epsilon(\bm e_{n \delta, \alpha_2}) \big).
\end{equation*}
Using the identity
\begin{equation*}
\log(1 - \zeta) = - \sum_{n = 1}^\infty \frac{\zeta^n}{n},
\end{equation*}
we get
\begin{equation*}
O_{11}(\zeta)^{-1} = \exp \big( \sum_{n = 1}^\infty \frac{\zeta^{n s}}{n} q^{n \bm K_1} \big).
\end{equation*}
We represent $\Phi(\zeta)$ as
\begin{equation*}
\Phi(\zeta) = \sum_{n = 1}^\infty \frac{\Phi_n}{n} \zeta^n
\end{equation*}
and see that equation (\ref{v11o}) is satisfied if
\begin{equation}
\Phi_n = q^{n \bm K_1} - \kappa_q n q^{-n} \epsilon(\bm e_{n \delta, \, \alpha_2}). \label{fn}
\end{equation}
It remains to demonstrate that equations (\ref{v22o}) and (\ref{v33o}) are also satisfied with this choice of $\Phi(\zeta)$.

It follows from (\ref{v11o}) that equation (\ref{v22o}) is equivalent to
\begin{equation*}
 V_{1 1}(\zeta)^{-1} V_{2 2}(\zeta) = O_{1 1}(\zeta)^{-1} O'_{2 2}(\zeta).
\end{equation*}
Using (\ref{v22}) and (\ref{v11}), we arrive at the equation
\begin{equation*}
\exp(-\kappa_q \bm \calE_{\delta, \, \alpha_1}(q^{-2} \zeta)) = O_{1 1}(\zeta)^{-1} O'_{2 2}(\zeta).
\end{equation*}
It follows from (\ref{edafda}) that
\begin{equation*}
\exp(-\kappa_q \bm \calE_{\delta, \, \alpha_1}(\zeta)) = 1 - \kappa_q \bm \calE'_{\delta, \, \alpha_1}(\zeta),
\end{equation*}
where
\begin{equation*}
\bm \calE'_{\delta, \, \alpha_1}(\zeta) = \sum_{n = 1}^\infty \epsilon(\bm e'_{n \delta, \, \alpha_1}) \zeta^n,
\end{equation*}
and we obtain the equation
\begin{equation}
\bm \calE'_{\delta, \, \alpha_1}(\zeta) = \frac{1}{\kappa_q} \big( 1 -  O_{1 1}(\zeta)^{-1} O'_{2 2}(\zeta)) \big) \label{epda1o}
\end{equation}
equivalent to equation (\ref{v22o}). It follows from (\ref{epnda}) that
\begin{equation}
\bm \calE'_{\delta, \, \alpha_1}(\zeta) = \zeta \big( \bm \calE_{\alpha_{1 2}}(\zeta) \epsilon(\bm e_{\delta - \alpha_{1 2}}) - q^2 \epsilon(\bm e_{\delta - \alpha_{1 2}}) \bm \calE_{\alpha_{1 2}}(\zeta) \big). \label{epda1}
\end{equation}
This relation uniquely determines the generating function $\bm \calE'_{\delta, \, \alpha_1}(\zeta)$. It satisfies equation (\ref{epda1}) by definition. One can verify that the right hand side of (\ref{epda1o}) also satisfies equation (\ref{epda1}). This means that equation (\ref{epda1o}) is valid. Therefore, equation (\ref{v22o}) is also valid.

It follows from (\ref{v22o}) that equation (\ref{v33o}) is equivalent to
\begin{equation*}
 V_{2 2}(\zeta)^{-1} V_{3 3}(\zeta) = O'_{2 2}(\zeta)^{-1} O''_{3 3}(\zeta).
\end{equation*}
Using (\ref{v33}) and (\ref{v22}) we come to the equation
\begin{equation*}
\exp(-\kappa_q \bm \calE_{\delta, \, \alpha_2}(q^{-2} \zeta)) = O'_{2 2}(\zeta)^{-1} O''_{3 3}(\zeta).
\end{equation*}
and, using (\ref{edafda}), obtain
\begin{equation*}
\exp(-\kappa_q \bm \calE_{\delta, \, \alpha_2}(\zeta)) = 1 - \kappa_q \bm \calE'_{\delta, \, \alpha_2}(\zeta),
\end{equation*}
where
\begin{equation*}
\bm \calE'_{\delta, \, \alpha_2}(\zeta) = \sum_{n = 1}^\infty \epsilon(\bm e'_{n \delta, \, \alpha_2}) \zeta^n.
\end{equation*}
Hence, we come to the equation
\begin{equation}
\bm \calE'_{\delta, \, \alpha_2}(\zeta) = \frac{1}{\kappa_q} \big( 1 -  O'_{2 2}(\zeta)^{-1} O''_{3 3}(\zeta)) \big) \label{epda2o}
\end{equation}
equivalent to equation (\ref{v33o}). It follows from (\ref{epnda}) that
\begin{equation}
\bm \calE'_{\delta, \, \alpha_2}(\zeta) = \zeta \big( \bm \calE_{\alpha_{2 3}}(\zeta) \epsilon(\bm e_{\delta - \alpha_{2 3}}) + \epsilon(\bm e_{\delta - \alpha_{2 3}}) \bm \calE_{\alpha_{2 3}}(\zeta) \big). \label{epda2}
\end{equation}
This relation uniquely determines the generating function $\bm \calE'_{\delta, \, \alpha_2}(\zeta)$. It satisfies equation (\ref{epda2}) by definition. One can verify that the right hand side of (\ref{epda2o}) also satisfies equation (\ref{epda2}) with the same initial condition. This means that equation (\ref{epda2o}) is valid. Therefore, equation (\ref{v33o}) is also valid.

\section{On central elements} \label{s:5}

It follows from the first relation of (\ref{edafda}) that
\begin{equation*}
\bm e_{n \delta, \, \alpha_i} = \sum_{k_1 + 2 k_2 + \cdots + k_n = n} \frac{\kappa_q^{k_1 + k_2 + \cdots + k_k - 1}(k_1 + k_2 + \cdots + k_n - 1)!}{k_1! k_2! \ldots k_n!} \, (\bm e'_{\delta, \, \alpha_i})^{k_1} \, (\bm e'_{2 \delta, \, \alpha_i})^{k_2} \, \ldots \, (\bm e'_{n \delta, \, \alpha_i})^{k_n}.
\end{equation*}
In particular,
\begin{align}
& \bm e_{\delta, \gamma} = \bm e'_{\delta, \gamma}, \label{edg1} \\
& \bm e_{2 \delta, \gamma} = \bm e'_{2 \delta, \gamma} + \kappa_q (\bm e'_{\delta, \gamma})^2 / 2, \\
& \bm e_{3 \delta, \gamma} = \bm e'_{3 \delta, \gamma} + \kappa_q \bm e'_{\delta, \gamma} \bm e'_{2 \delta, \gamma} + \kappa_q^2  (\bm e'_{\delta, \gamma})^3 / 3, \\
& \bm e_{4 \delta, \gamma} = \bm e'_{4 \delta, \gamma} + \kappa_q \bm e'_{\delta, \gamma} \bm e'_{3 \delta, \gamma} + \kappa_q (\bm e'_{2 \delta, \gamma})^2 + \kappa_q^2  (\bm e'_{\delta, \gamma})^2 \bm e'_{2 \delta, \gamma} + \kappa_q^3  (\bm e'_{\delta, \gamma})^4 / 4. \label{edg4}
\end{align}

We denote
\begin{equation*}
\obmE_{1 2} = \bm E_{1 2}, \qquad \obmE_{1 3} = \bm E_1 \bm E_2 - q^{-1} \bm E_2 \bm E_1, \qquad \obmE_{2 3} = \bm E_{2 3}.
\end{equation*}
Using the first equations of (\ref{epda}) and (\ref{epnda}), we obtain
\begin{align*}
\epsilon(\bm e'_{\delta, \, \alpha_2}) & = -\kappa_q^{-1} (q^{-1} q^{2 \bm K_2} - q^{-1} q^{-2 \bm K_3}) - \kappa_q (\bm F_{1 2} \, \obmE_{1 2} \, q^{\bm K_1 + \bm K_2} + q^2 \bm F_{1 3} \, \obmE_{1 3} \, q^{\bm K_1 - \bm K_3} + \bm F_{2 3} \, \obmE_{2 3} \, q^{\bm K_2 - \bm K_3}), \\
\epsilon(\bm e'_{2 \delta, \, \alpha_2}) & = - \kappa_q^{-1} (q^{-2} q^{4 \bm K_2} - q^{-2} q^{2 \bm K_2 - 2 \bm K_3}) - \kappa_q ([2]_q q^{_2} \bm F_{1 2} \, \obmE_{1 2} \, q^{\bm K_1 + 2 \bm K_2} + q^3 \bm F_{1 2} \, \obmE_{1 2} \, q^{3 \bm K_1 + \bm K_2} \\
& - q^{-1} \bm F_{1 2} \, \obmE_{1 2} \, q^{\bm K_1 + \bm K_2 - 2 \bm K_3}+ q^5 \bm F_{1 3} \, \obmE_{1 3} \, q^{3 \bm K_1 - \bm K_3} + q \bm F_{1 3} \, \obmE_{1 3} \, q^{\bm K_1 + 2 \bm K_2 - 1 \bm K_3} \\
& + [2]_q \bm F_{2 3} \, \obmE_{2 3} \, q^{3 \bm K_2 -\bm K_3} ) - \kappa_q^2 (q^{-1} \bm F_{1 2} \, \bm F_{2 3} \, \obmE_{1 3} \, q^{\bm K_1 + 2 \bm K_2 - \bm K_3} + q^2 \bm F_{1 3} \obmE_{1 2} \, \obmE_{1 3} \, q^{\bm K_1 + 2 \bm K_3 - \bm K_3}) \\
& + \kappa_q^3 ([2]_q q^3 \bm F_{1 2} \, \bm F_{1 3} \, \obmE_{1 2} \, \obmE_{1 3} \, q^{2 \bm K_1 + \bm K_2 - \bm K_3} + \bm F_{1 2} \, \bm F_{2 3} \, \obmE_{1 2}  \, \obmE_{2 3} q^{2 \bm K_1 + \bm K_2 - \bm K_3} + (F_{1 2})^2 (\obmE_{1 2})^2 q^{2 \bm K_1 + 2 \bm K_2})
\end{align*}
and the relation (\ref{fn}) gives
\begin{align*}
\Phi_1 = {} & q^{2 \bm K_1} + q^{-2} q^{2 \bm K_2} - q^{-2} q^{-2 \bm K_3} + \kappa_q^2 (q^{-1} \bm F_{1 2} \, \obmE_{1 2} \, q^{\bm K_1 + \bm K_2} + q \bm F_{1 3} \, \obmE_{1 3} q^{\bm K_1 - \bm K_3} \\
& + q^{-1} \bm \, F_{2 3} \, \obmE_{2 3} \, q^{\bm K_2 - \bm K_3}), \\
\Phi_2 = {} & q^{4 \bm K_1} + q^{-4} q^{4 \bm K_2} - q^{-4} q^{-4 \bm K_3} + \kappa_q^2 [2]_q (q^{-4} \bm F_{1 2} \, \obmE_{1 2} \, q^{\bm K_1 + 3 \bm K_2} + \bm F_{1 2} \, \obmE_{1 2} \, q^{\bm 3 K_1 + \bm K_2} \\
& + \bm F_{1 3} \, \obmE_{1 3} \, q^{\bm K_1 - 3 \bm K_3} + q^2 \bm F_{1 3} \, \obmE_{1 3} \, q^{\bm 3 K_1 - \bm K_3} + q^{-2} \bm F_{2 3} \, \obmE_{2 3} \, q^{\bm K_2 - 3 \bm K_3} + q^{-2} \bm F_{2 3} \, \obmE_{2 3} \, q^{3 \bm K_2 - \bm K_3}) \\
& + \kappa_q^3 [2]_q (q^{-2} \bm F_{1 2} \, \bm F_{2 3} \, \obmE_{1 3} \, q^{\bm K_1 + 2 \bm K_2 - \bm K_3} + q^{-1} \bm F_{1 3} \, \obmE_{1 2} \, \obmE_{2 3} \, q^{\bm K_1 + 2 \bm K_2 - \bm K_3}) \\
& + \kappa_q^4 [2]_q (q \bm F_{1 2} \, \bm F_{1 3} \, \obmE_{1 2} \, \obmE_{1 3} \, q^{2 \bm K_1 + \bm K_2 - \bm K_3} + q \bm F_{1 3} \, \bm F_{2 3} \, \obmE_{1 3} \, \obmE_{2 3} \, q^{\bm K_1 + \bm K_2 - 2 \bm K_3}) \\
& + \kappa_q^4 \, q^{-2} (\bm F_{1 2})^2 (\obmE_{1 2})^2 q^{2 \bm K_1 + 2 \bm K_2}.
\end{align*}
One can demonstrate that the elements $\Phi_1$ and $\Phi_2$ belong to the center of $\uqglto$.

It is evident that the eigenvalues of $F_1$ and $F_2$ in the irreducible finite dimensional representation $V^\lambda$ of highest weight $\lambda = (\lambda_1, \lambda_2, \lambda_3)$ are
\begin{equation*}
\chi_\lambda(\Phi_1) = q^{2 \bm \lambda_1} + q^{2 \bm \lambda_2 - 2} - q^{-2 \bm \lambda_3 - 2}, \qquad \chi_\lambda(\Phi_2) = q^{4 \bm \lambda_1} + q^{4 \bm \lambda_2 - 4} - q^{-4 \bm \lambda_3 - 4}.
\end{equation*}
It is natural to assume that for a general $n$ we have
\begin{equation*}
\chi_\lambda(\Phi_n) = q^{2 n \bm \lambda_1} + q^{2 n (\bm \lambda_2 - 1)} - q^{-2 n (\bm \lambda_3 + 1)}.
\end{equation*}
We checked the validity of this conjecture also for $F_3$ and $F_4$.

In the papers \cite{ZhaGou91, Zha91, Zha92}, the method of partial trace proposed by Drinfeld \cite{Dri90} was used to find the central elements of quantum superalgebras $\uqglmn$.\footnote{The paper \cite{Dri90} addresses the case of ordinary quantum algebras. The generalization to the quantum superalgebras is straightforward. The Russian original of the paper \cite{Dri90} is: {\Russian В. Г. Дринфельд, \emph{О почти кокоммутативных алгебрах Хопфа}, Алгебра и анализ}, \textbf{1:2} (1989), 30--46.} The essence of this approach is as follows.

Let $A$ be an almost cocommutative Hopf superalgebra with the comultiplication $\Delta$ and the antipode $S$. It means that there exists an invertible even element $\calR \in A \otimes A$ such that
\begin{equation*}
(\sigma \circ \Delta)(a) = \calR \, \Delta(a) \calR^{-1} \label{sda}
\end{equation*}
for any $a \in A$. Here $\sigma$ is an automorphism of $A \otimes A$ defined on the homogeneous elements as
\begin{equation*}
\sigma(a_1 \otimes a_2) = (-1)^{[a_1] [a_2]} a_2 \otimes a_1.
\end{equation*}
Let $B$ be an even element of $A \otimes A$ such that
\begin{equation}
[\Delta(a), B] = 0 \label{dab}
\end{equation}
for all $a \in A$. If an element $\lambda \in A^*$ satisfies the equation
\begin{equation}
\lambda(a b) = (-1)^{[a] [b]} \lambda(b S^2(a)) \label{lab}
\end{equation}
for any $a, b \in A$, then $(1 \otimes \lambda)(B)$ belongs to the center of $A$.

The quantum superalgebra $\uqglto$ is quasitriangular. Hence, it is almost cocommutative. We can demonstrate that for
\begin{equation*}
\bm K = 2 \bm K_2 + 2 \bm K_3
\end{equation*}
we have
\begin{equation*}
S^2(a) = q^{\bm K} a q^{-\bm K}
\end{equation*}
for any $a \in \uqglto$. Using this equation, we rewrite equation (\ref{lab}) as
\begin{equation*}
\lambda(a b q^{-\bm K}) = (-1)^{[a] [b]} \lambda(b a q^{-\bm K}). 
\end{equation*}
It is clear that for the element $\mu \in \uqglto^*$ defined by the equation
\begin{equation*}
\mu(a) = \lambda(a q^{-\bm K})
\end{equation*}
we have
\begin{equation}
\mu(a b) = (-1)^{[a] [b]} \mu(b a). \label{mab}
\end{equation}
Now, it is not difficult to understand that, for a given element $B \in \uqglto \otimes \uqglto$ satisfying equation (\ref{dab}) and any element $\mu \in \uqglto^*$ satisfying equation (\ref{mab}), the element
\begin{equation*}
C = (1 \otimes \mu)(B (1 \otimes q^{\bm K})
\end{equation*}
belongs to the center of $\uqglto$. It is evident that if $B$ satisfies equation (\ref{dab}), then $B^n$ for any integer $n > 1$ satisfies the same equation, and
\begin{equation*}
C_n = (1 \otimes \mu)(B^n (1 \otimes q^{\bm K})).
\end{equation*}
belongs to the center of $\uqglto$.

It is easy to get convinced that $B = \sigma(\calR) \calR$ satisfies equation (\ref{dab}), and $\mu = \str \circ \pi$ satisfies equation (\ref{mab}). Hence, the elements
\begin{equation*}
C_n = (1 \otimes (\str \circ \pi))((\sigma(\calR) \calR)^n (1 \otimes q^{\bm K}))
\end{equation*} 
belong to the center of $\uqglto$.

The universal $R$-matrix for the quantum superalgebra $\uqglto$ has the form \cite{KhoTol91}
\begin{equation*}
\calR = \exp_{q^2}(- \kappa_q \, \bm E_{1 2} \otimes \bm E_{2 1}) \exp(\kappa_q \, \bm E_{1 3} \otimes \bm E_{3 1}) \exp(\kappa_q \, \bm E_{2 3} \otimes \bm E_{3 2}) \exp(-\sum_{i = 1}^3 (-1)^{[i]} \bm K_i \otimes \bm K_j),
\end{equation*}
and we obtain
\begin{equation*}
C_n = (1 \otimes \str)((M^\sigma M)^n (1 \otimes\pi(q^{\bm K}))),
\end{equation*}
where
\begin{multline*}
M = (1 \otimes \pi)(\calR) = (1 \otimes \mathbbm 1 - \kappa_q \, \bm E_{1 2} \otimes \bbM_{2 1} + \kappa_q \, \bm E_{1 3} \otimes \bbM_{3 1} + \kappa_q \, \bm E_{2 3} \otimes \bbM_{3 2}) \\
\times (q^{- \bm K_1} \otimes \bbM_{1 1} + q^{- \bm K_2} \otimes \bbM_{2 2} + q^{\bm K_3} \otimes \bbM_{3 3}),
\end{multline*}
\begin{multline*}
M^\sigma = (1 \otimes \pi)(\sigma(\calR)) = (1 \otimes \mathbbm 1 - \kappa_q \, \bm E_{2 1} \otimes \bbM_{1 2} + \kappa_q \, \bm E_{3 1} \otimes \bbM_{1 3} + \kappa_q \, \bm E_{3 2} \otimes \bbM_{3 2}) \\
\times (q^{- \bm K_1} \otimes \bbM_{1 1} + q^{- \bm K_2} \otimes \bbM_{2 2} + q^{\bm K_3} \otimes \bbM_{3 3}),
\end{multline*}
and
\begin{equation*}
\pi(q^{\bm K}) = \bbM_{1 1} + q^2 \bbM_{2 2} + q^2 \bbM_{3 3}.
\end{equation*}

We define the Cartan anti-involution $\omega$ for $\uqglto$ by the relations
\begin{equation*}
\omega(\bm E_i) = \bm F_i, \qquad \omega(\bm F_i) = \bm E_i, \qquad \omega(\bm H_i) = \bm H_i
\end{equation*}
together with the rule $\omega(\hbar) = -\hbar$ implying that $\omega(q) = q^{-1}$. By definition, for any integer $n \ge 1$, the element
\begin{equation*}
\widetilde C_n = \omega(C_n)
\end{equation*}
belongs to the center of $\uqglto$. It seems plausible that the elements of $F_n$ are algebraically dependent on the elements $\widetilde C_n$. In particular, we have
\begin{align*}
& \Phi_1 = \widetilde C_1, \\
& \Phi_2 = \frac{1}{[2]_q} (2 q \widetilde C_2 - \kappa_q (\widetilde C_1)^2), \\
& \Phi_3 = \frac{1}{[3]_q} (3 q^2  \widetilde C_3 - 3 \kappa_q q \widetilde C_2 \widetilde C_1 + \kappa_q^2 (\widetilde C_1)^3), \\
& \Phi_4 = \frac{1}{[4]_q} (4 q^3 \widetilde C_4 - 4 \kappa_q q^2 \widetilde C_3 \widetilde C_1 - 2 \kappa_q q^2 \widetilde C_2 \widetilde C_2 + 4 \kappa_q^2 q \widetilde C_2 (\widetilde C_1)^2 - \kappa_q^3 (\widetilde C_1)^4).
\end{align*}
It is instructive to compare these relations with the relations (\ref{edg1})--(\ref{edg4}).

\section{Conclusions} \label{s:6}

We derived the explicit expression for the basic monodromy operator for the quantum loop superalgebra $\uqlslto$. Two significant additional results emerged from this derivation:
simple expressions for the generating functions of the the images of the root vectors of $\uqlslto$ under the Drinfeld homomorphism and explicit expressions for certain central elements of the quantum superalgebra $\uqglto$. We established the relation between the derived central elements and those obtained using the Drinfeld partial trace method. An advantage of the derived central elements is the simple form of their eigenvalues. The expressions for the generating functions of the root vectors could be used to find the $\ell$-weights of the evaluation representations of the quantum loop superalgebra $\uqlslto$.

\appendix

\section{Superalgebras} \label{a:a}

Let $V$ be a $\bbZ_2$-graded vector space so that
\begin{equation*}
V = V_\oz \oplus V_\oo.
\end{equation*}
The elements of the subspaces $V_\oz$ and $V_{\overline 1}$ are said to be even and odd, respectively. An element belonging to $V_\oz$ or $V_\oo$ is said to be homogeneous. The parity of a homogeneous element $v$, denoted as $[v]$, is $\overline 0$ or $\overline 1$ according to whether it is in $V_\oz$ or $V_\oo$.\footnote{Throughout the paper, we use the following convenient convention. If $v$ is an element of a $\bbZ_2$-graded vector space and $[v]$ appears in some formula or expression, then $v$ is assumed to be homogeneous.} For a finite dimensional $\bbZ_2$-graded vector space $V$, the pair $(\dim V_\oz, \, \dim V_\oo)$ is called the dimension of $V$. The tensor product of two $\bbZ_2$-graded vector spaces $V$ and $W$ is considered as a $\bbZ_2$-graded vector space with
\begin{equation*}
(V \otimes W)_\alpha = \bigoplus_{\beta + \gamma = \alpha} V_\beta \otimes W_\gamma.
\end{equation*}
It is clear that
\begin{equation*}
[v \otimes w] = [v] + [w].
\end{equation*}

An algebra $A$ is called a superalgebra if it is $\bbZ_2$-graded as a vector space and
\begin{equation*}
A_\alpha A_\beta \subset A_{\alpha + \beta}, \qquad \alpha, \beta \in \bbZ_2.
\end{equation*}
It follows that for any two homogeneous elements $a_1, a_2 \in A$, the product $a_1 a_2$ is homogeneous and
\begin{equation*}
[a_1 a_2] = [a_1] + [a_2].
\end{equation*}
The tensor product of two associative superalgebras $A$ and $B$ is also an associative superalgebra. Its multiplication is defined by the equation
\begin{equation*}
(a_1 \otimes b_1)(a_2 \otimes b_2) = (-1)^{[b_1][a_2]} (a_1 a_2) \otimes (b_1 b_2)
\end{equation*}
for any homogeneous elements $a_1, a_2 \in A$ and $b_1, b_2 \in B$. The multiplication in the tensor product of any finite number of associative superalgebras is defined recursively.

Let $V$ be a $\bbZ_2$-graded vector space and $M \in \End(V)$. We have
\begin{equation*}
M v = M v_\oz + M v_\oo = (M v_\oz)_\oz + (M v _\oz)_\oo + (M v_\oo)_\oz +(M v_\oo)_\oo.
\end{equation*}
We define
\begin{equation*}
M_\oz v = (M v_\oz)_\oz + (M v_\oo)_\oo, \qquad M_\oo v = (M v _\oz)_\oo + (M v_\oo)_\oz
\end{equation*}
so that
\begin{equation*}
M = M_\oz + M_\oo.
\end{equation*}
The above decomposition is unique. Therefore, the associative algebra $\End(V)$ has a natural structure of a $\bbZ_2$-graded vector space. It is clear that
\begin{equation*}
[M v] = [M] + [v].
\end{equation*}
Furthermore, it can be shown that for any two homogeneous elements $M, N \in \End(V)$ we have
\begin{equation*}
[M N] = [M] + [N].
\end{equation*}
This implies that $\End(V)$ is an associative superalgebra.

A Lie superalgebra $\gothg$ is a superalgebra with an operation $[\, \cdot \, , \, \cdot \,]$ satisfying the conditions
\begin{gather*}
[a, \, b] = (-1)^{[a] [b]} [b, \, a], \\
(-1)^{[a] [c]}[a, \, [b, \, c]] + (-1)^{[b] [a]}[b, \, [c, \, a]] + (-1)^{[c] [b]}[c, \, [a, \, b]] = 0.
\end{gather*}
for all homogeneous elements $a, b, c$ in $\gothg$. Here and throughout the paper, it is assumed that
\begin{equation*}
(-1)^\oz =  1, \qquad (-1)^\oo =  - 1.
\end{equation*}
It is clear that
\begin{equation*}
[\gothg_\oz, \, \gothg_\oz] \subset \gothg_\oz, \qquad [\gothg_\oz, \, \gothg_\oo] \subset \gothg_\oo.
\end{equation*}
It means that $\gothg_\oz$ is a Lie algebra, and $\gothg_\oo$ is a $\gothg_\oz$-module.

Given an associative superalgebra $A$, we define the supercommutator on homogeneous elements of $A$ by
\begin{equation*}
[a_1, \, a_2] =  a_1 a_2 - (-1)^{[a_1] [a_2]} a_2 a_1
\end{equation*}
and then extend the definition by linearity to all elements of $A$. With respect to the supercommutator, the algebra $A$ is a Lie superalgebra. In this way, for any $\bbZ_2$-graded vector space $V$, starting with the superalgebra $\End(V)$, we obtain a Lie superalgebra called the general linear superalgebra and denoted by $\mathfrak{gl}(V)$.

\end{document}